\documentclass[11pt, a4paper]{article}

\usepackage{amsmath,amssymb, amsthm}
\usepackage{mathtools}
\usepackage[utf8]{inputenc}
\usepackage{mathrsfs}

\usepackage{fullpage}

\usepackage{authblk}
    
\usepackage{dsfont}
\usepackage{bbm}
    
\usepackage{graphicx}
\graphicspath{ {images/} }

\usepackage{wrapfig}

\usepackage{verbatim}

\usepackage{xcolor}

\usepackage{hyperref}

\hypersetup{
	colorlinks   = true, 
	urlcolor     = {blue!90!black}, 
	linkcolor    = {blue!90!black}, 
	citecolor   = {red!90!black} 
}

\usepackage{xpatch}
\xpatchcmd{\proof}{\itshape}{\normalfont\proofnamefont}{}{}

\newcommand{\proofnamefont}{\bfseries}

\newcommand{\R}{\mathbb{R}}
\newcommand{\N}{\mathbb{N}}
\newcommand{\Z}{\mathbb{Z}}
\renewcommand{\d}{{\mathrm{d}}} 

\newcommand{\m}{\mathbf{m}}
\renewcommand{\P}{\mathbb{P}}
\newcommand{\E}{\mathbb{E}}
\newcommand{\PP}[1]{\mathbb{P}\left\{#1\right\}}

\newcommand{\I}{{\rm{I}}}
\newcommand{\W}{{\rm{W}}}
\newcommand{\0}{\mathbf{0}}
\newcommand{\1}{\mathds{1}}
\newcommand{\tip}[1]{\mathrm{tip}\left(#1\right)}
\newcommand{\A}{\mathcal{A} }
\newcommand{\rmu}{{ \rm{u} }  }

\theoremstyle{plain}
\newtheorem{thm}{Theorem}[section]
\newtheorem{prop}[thm]{Proposition}
\newtheorem{lem}[thm]{Lemma}
\newtheorem{cor}[thm]{Corollary}
\newtheorem{rmk}[thm]{Remark}

\theoremstyle{definition}

\newtheorem{defi}[thm]{Definition}

\author[1,2]{Viktor Bezborodov \thanks{Email: \texttt{viktor.bezborodov@pwr.edu.pl}}} 
\author[2]{
 Tyll Krueger \thanks{Email: \texttt{tyll.krueger@pwr.wroc.pl}}}

\affil[1]{
	{University of Goettingen,  Institute for Mathematical Stochastics
}}
\affil[2]{
	{Wroc\l{}aw University of Science and Technology
	 }}

\title{Linear and superlinear spread for 
stochastic combustion growth process}

\begin{document}

\maketitle

\begin{abstract}
 
Consider a stochastic growth model on $\Z ^\d$.
Start with some active particle at the origin and
sleeping particles
elsewhere.
The initial number of particles at  $x \in \Z ^\d$
is
$\eta(x)$, where $\eta (x)$
are independent random variables distributed according to $\mu$.
Active particles perform a simple continuous-time 
random walk
while sleeping particles
stay put until 
the first arrival of an active particle to their location.
Upon the arrival all sleeping particles at the site
activate at once and start moving
according to their own simple random walks.
The aim of this paper is to give
conditions on $\mu$
under which the spread of the process
is linear or 
faster than 
linear.  
The proofs 
rely on
comparison to various percolation models.
 
\end{abstract}



%
%

%

\section{Introduction}

At time $t = 0$ there are  $\eta (x)$  particles at $x\in \Z ^\d$,
where the random variables $\{\eta(x) \}_{x \in \Z ^\d}$ are independent and identically distributed
according to a distribution $\mu$ on $\Z_+ := \N \cup \{0\}$, 
$\N = \{1,2,3,...\}$,
and $\d \in \N$ is the dimension.
The particles at the origin are active while all other particles are dormant (sleeping).
Active particles perform a simple continuous-time random walk
independently of all other particles.
Sleeping particles stay still until the first
arrival of an active particle to their location;
upon arrival they become active and start their own simple random walks.
We exclude a trivial case
and  assume throughout that $\mu(0) < 1$.
Denote by $\A_t$  the set of sites visited by 
an active particle by the time $t$.
In this paper we 
investigate the various conditions on $\mu$
ensuring that the system spreads linearly with time,
or that the system spreads faster than linearly with time.

\begin{defi}\label{def linear spread}
	We say that the spread is linear,
	or the spread rate is linear,
	or the system spreads linearly,
	if there exists a constant $C > 0$
	such that a.s.
	\begin{equation} \label{flick}
	\A _t \subset C [-t,t] ^\d
	\  \  \  \text{for large } t > 0.
	\end{equation}
	
	If a.s. \eqref{flick} does not hold for any $C >0$,
	the spread (rate) is said to be superlinear, or faster than linear;
	the system spreads faster than linearly with time. 
	In other words, the spread is superlinear if 
	for every $C > 0$ the set 
	$
	\{ t \geq 0: \A _t \not\subset C [-t,t] ^\d \}
	$ 
	is unbounded a.s.
\end{defi}

It can be shown using a $0-1$
type of argument that
in the one-dimensional case 
the spread  is superlinear 
if and only if $\limsup\limits_{ t \to \infty }\frac{\sup \A_t}{t} = \infty$
(we do not prove it in this paper because it is not needed in our proof of 
Theorem \ref{thm main}).
Sometimes instead of the linear spread
the phrase `the linear growth' (or `the linear growth rate')
is used to describe \eqref{flick}.
Let $\0$ be the origin of $\Z ^\d$.
In the case $\eta(\0)=0$ a single active particle
is added to the origin to prevent 
a possible absence of active particles.
We collect  our principal results in the next theorem.

\begin{thm} \label{thm main} 
	Consider the stochastic combustion growth process. 
	\begin{itemize}
		\item [](i)  Assume that
		for some $B > 1$,
		\begin{equation}\label{wack}
		\sum\limits _{m \in \N} \left[
		\mu \big( [ B^{m} , \infty ) \big) 
		\right] ^{\frac 1 \d}  < \infty.
		\end{equation}
		Then  the spread is linear.

		\item [](ii)
		Assume that for every $B >1$,
		\begin{equation}  \label{obtuse}
		\sum\limits _{m=1} ^\infty \prod\limits _{n = 1} ^ m 
		\mu \left( [0, B ^ n] \right) < \infty.
		\end{equation}
		Then  the spread is superlinear.

		\item  [] (iii)  Assume that  for some $B>1$,
		\begin{equation}\label{enclave}
		\sum\limits _{n \in \N}   \mu \big( [ B ^ {n \ln ^2 n}, \infty) \big)  = \infty.
		\end{equation}
		
		Then  the spread is superlinear.

	\end{itemize}
	
\end{thm}

We postpone a discussion
on the cases
not covered 
by Theorem \ref{thm main}
to Remark \ref{snorkel}.
The asymptotic shape
of the discrete time version of 
the model is given in \cite{shapeFrog, shapeFrogRandom}.
The 
stochastic combustion growth process
was introduced in
\cite{stocCombust},
and it is also known
as a  model of $X + Y \to 2X$ \cite{CQR09, BR10}.
In \cite{stocCombust}
$\mu $ is  the delta measure at $1$;
that is, at the beginning
there is exactly one particle per site. 
A shape theorem
is proven
in that paper,
and it is shown that 
the distribution of the number of
particles in visited sites
converges to the product Poisson measure
with parameter 1 (see \cite{stocCombust} for the precise formulation).
For the one-dimensional 
stochastic combustion growth process
further extensions have been obtained.
A central limit theorem for the front of the 
process 
with a fixed number of sleeping particles per site
is given in \cite{CQR09}.
For a slightly modified
model  in \cite{CQR07}
a shape theorem and a central limit theorem 
for the front are established. 
In \cite{Gia09} a lower bound on the 
speed of the front in  a one-dimensional 
model with Poisson initial condition is given.

Extensions 
with sleeping particles replaced by moving particles of different type
are treated in \cite{KS05, KS08, BR16, BT20}.
In particular, in \cite{KS05} 
the linear growth is established, while in \cite{KS08}
a shape theorem   is proven. 
Further discussion
of this and related models 
takes place 
in \cite{KRS12}.
In \cite{BR16} 
a central limit theorem is obtained
for the 
front in the one-dimensional models
with mobile particles of two types.
Linear spread of an infection spread
on a zero range process is proven in \cite{BT20}.

Perhaps the fist particle growth model
where the spread rate 
can be linear or superlinear depending
on  parameters 
was the branching random walk. 
The spread rate is linear for the branching random 
walk with exponential tails \cite{Dur79}
and exponential for the branching random walk
with polynomial tails \cite{Dur83}. 
An intermediate case 
can result in a polynomial spread rate
\cite{Gan00}.
The
exact expression of the 
speed of a
discrete-space
branching random walk with an exponential moment condition 
can be found in \cite{Big95}.
The model in  \cite{trunc_and_crop} 
is
the branching random walk
with the additional restriction 
that the birth rate at any spatial
location 
cannot exceed one. It is shown 
that this model spreads linearly
provided the tails 
of the dispersion kernel are lighter than $ |x|^{-4 - \varepsilon}$
for some $\varepsilon >0$.

In this paper we establish conditions 
for the linear or superlinear spread
using comparison to certain
percolation 
and percolation-like models. 
Comparing a growth process with a percolation model
is not uncommon. 
Famously, the renormalization
procedure
in the proof of the shape theorem
for the contact process 
is a stepping stone 
toward comparison to  oriented percolation  \cite{DG82, Dur91}.
The comparison of the contact process to a simpler growth model
via the renormalization
proves useful in various contexts \cite[Chapter 2]{Liggettbook2}.
Renormalization procedure
and comparison to another process
is  used in 
percolation theory itself \cite{ADH15}.

Domination by other models
is a common technique when 
proving a linear growth
of a certain stochastic process.
In \cite{Sidoravicius19} 
a certain aggregation process
is shown to grow linearly with time
via comparison with a two type
first passage percolation process. 
In \cite{GM08}
and in this paper 
the growth process
is compared to a greedy paths model \cite{CGGK93}
(more precisely, in \cite{GM08}
a continuous-space equivalent is used).
The linear speed of the continuous-space growth models
can often be deduced from the linear speed
of similar discrete-space models
\cite{Dei03}, \cite{shapenodeath}.
In \cite{Brownian_Frog_18}
the Brownian frog model is 
dominated by a certain specially designed
branching process. 
In \cite{trunc_and_crop} 
the system viewed from its tip
is dominated by another more amenable
to analysis process.

The discrete-time version
of the model described at the beginning
is known as the frog model.
Respectively,
the stochastic combustion growth process
can be seen as the continuous-time frog model. 
The discrete-time model has been 
an active research subject in recent years.
In discrete time the  model 
cannot grow faster than linearly
as the set of  visited sites
is always contained in  $t\mathcal{D}$, 
where $\mathcal{D} = \{ (x_1, \dots, x_\d) : |x_1| + \dots + |x_\d| \leq 1 \}$
and $t = 0,1,\dots$.
The shape theorem was proven 
in \cite{shapeFrog} and \cite{shapeFrogRandom}.
The question whether $\mathcal{D}$
can be a limiting shape 
for distributions $\mu$
with sufficiently heavy
tails was answered positively 
in \cite{shapeFrogRandom}.
Recent papers \cite{DHL19} and \cite{benjamini2020}
provide an overview of 
other research on this model.  
The transitivity and recurrence properties
of the frog model
attract considerable attention \cite{RecFrog18,   RecFrog17, HJJ16, KZ17, GNR17, GTW22}.
In \cite{Zer18} 
a recurrence criterion is obtained for an asymmetric frog model
with particles removed after a geometric number of steps. 
For various  finite graphs
the asymptotics of the first moment
when all sites are visited
is established in 
\cite{benjamini2020}.
Continuity of the asymptotic
shape of the discrete-time frog model
with respect to the measure $\mu$
is established in  \cite{Kub20}.
The variance of the passage times is sublinear
\cite{CN19}, which implies
that a central limit theorem does not hold.
In \cite{DHL19} a possibility of  co-existence
in a two-type frog model with  lazy walkers was demonstrated. 
The co-existence result was
extended in dimension one
to walkers with different step probabilities in \cite{HK22}.

In the Brownian frog model
the particles perform a Brownian
motion instead of a simple random walk.
Naturally, the process evolves in continuous time. 
A shape theorem and an asymptotic density results
were obtained in \cite{Brownian_Frog_18},
while conditions for transience of the one-dimensional 
version
are established in \cite{Brownian_Frog_17} (active particles
in \cite{Brownian_Frog_17}  have a leftward drift, 
thus it is possible that all of them escape to $-\infty$ and 
the origin is
not  being
visited starting from some positive time onwards).

The paper is organized as follows.
Further definitions, notation, comments, 
the structural description of the paper,
and  the proof ideas
are collected in Section \ref{sec 2}.
In Section \ref{sec TADBP percolation}
the properties of an auxiliary percolation model 
are given.
The proof of Theorem \ref{thm main}
is spread across Sections \ref{sec proofs directly from TADBP},
\ref{sec proofs TADBP indirect},
and  \ref{sec proofs animals}.
In Section \ref{sec convergence properties}
the
independence of convergence properties 
of series in (ii)
and (iii) of Theorem \ref{thm main} is discussed.

\section{Further discussion and some ideas of the proof}\label{sec 2}

Let us first discuss the proof structure.
The proof of (i) of Theorem 
\ref{thm main} in dimension $\d = 1$
follows virtually the same steps
as the proof  
of (ii) of Theorem 
\ref{thm main}, 
whereas in 
case $\d \geq 2$
the proof of 
(i) is very different.
To streamline the following discussions
and the proofs we now formulate two
theorems which are exactly 
(i) of Theorem 
\ref{thm main}
in cases $\d = 1$
and $\d \geq 2$.
Thus, the first theorem of this section gives sufficient conditions for the linear spread 
of the one-dimensional system. 

\begin{thm} \label{thm linear growth d = 1}
	Let  $\d = 1$ 
	and
	assume 
	\begin{equation} \label{egregious}
	\sum\limits _{k = 1} ^\infty  \mu(k) \ln k  < \infty. 
	\end{equation}
	Then 	the spread is linear in time.
\end{thm}

We note that for $\d = 1$, \eqref{egregious} 
is equivalent to \eqref{wack}.
The next theorem gives sufficient conditions
for the linear spread in dimensions $\d \geq 2$.
As expected,
the assumptions are stronger than in the case $\d = 1$.

\begin{thm}\label{thm linear growth d geq 2}
	Let $\d \geq  2$ and assume that
	for every $B > 1$,
	\begin{equation} \label{ad nauseam}
	\sum\limits _{m \in \N} \left[
	\mu(  [ B^{m} , \infty ) )
	\right] ^{\frac 1 \d}  < \infty.
	\end{equation}
	Then  the spread is linear.
\end{thm}

When combined,
the statements of Theorems \ref{thm linear growth d = 1}
and \ref{thm linear growth d geq 2}
exactly make up (i)
of Theorem \ref{thm main}.
The proof of Theorem \ref{thm linear growth d = 1}
can be found in Section \ref{sec proofs directly from TADBP}.
The proofs of
(ii) and (iii)
of
Theorem
\ref{thm main}
are located in Sections \ref{sec proofs directly from TADBP}
and \ref{sec proofs TADBP indirect}, respectively. 
A brief discussion on the cases not covered by Theorem \ref{thm main}
can be found at the end of Section \ref{sec proofs TADBP indirect} in Remark \ref{snorkel}.
The proof of Theorem \ref{thm linear growth d geq 2}
is contained in Section \ref{sec proofs animals}.

The following proposition 
is useful because 
it shows that
it is enough
to prove 
that the spread 
is
superlinear 
for $\d = 1$ only.
In particular, in the proofs of
(ii) and (iii)
of
Theorem
\ref{thm main}
it is sufficient to consider the dimension $\d = 1$ only.
A monotonicity in dimension 
of this kind
appears in \cite{stocCombust}.

\begin{prop} \label{prop superlinear growth higher dimension}
	Assume  $\mu$ is such that the spread of 
	the corresponding one-dimensional
	model is superlinear a.s.
	Then a.s. the spread 
	is superlinear for $\d \geq 2$ as well.
\end{prop}

\begin{rmk}
	Proposition \ref{prop superlinear growth higher dimension}
	can be extended to dimensions $d_1$ and $d_2$
	with
	$1 \leq \d_1 < \d_2$.
\end{rmk}

The proof of Proposition \ref{prop superlinear growth higher dimension} can be found on Page \pageref{page proof prop superlinear growth higher dimension}.
Now we formulate a shape theorem,
which in this case is a consequence to linear growth.
For two sets $\mathbb{A}, \mathbb{B}$, let their 
sum be defined in the usual way $\mathbb{A} + \mathbb{B} = \{ a + b: a \in \mathbb{A}, b 
\in \mathbb{B} \}$.

\begin{cor}[{Shape Theorem}]
	Assume $\d$ and $\mu$ satisfy conditions of Theorem \ref{thm linear growth d = 1}
	or Theorem \ref{thm linear growth d geq 2}.	
	There exists a bounded non-empty convex set
	$\boldsymbol{\rm{A}}$ such that 
	for any $\varepsilon \in (0,1)$, 
	\begin{equation}
	(1-\varepsilon) \boldsymbol{\rm{A}} 
	\subset 
	\frac{\A_t + [-\frac 12, \frac 12] ^\d}{t}
	\subset 
	(1+\varepsilon) \boldsymbol{\rm{A}} 
	\end{equation}
	for all sufficiently large $t$. 
\end{cor}

We do not prove the shape theorem in this paper and refer instead to Section 3 
of \cite{shapeFrogRandom}.
The authors of that paper point out that the shape theorem 
for the discrete-time frog model proven 
in that paper holds for the continuous-time
version too, provided that 
the faster than linear spread is ruled out
-- and this is exactly 
the conclusions of Theorem \ref{thm linear growth d = 1}
and Theorem \ref{thm linear growth d geq 2}.

As mentioned in the introduction,
we make use of   auxiliary
models. 
In the proofs of  Theorem   \ref{thm linear growth d = 1},
and (ii) and (iii)
of
Theorem
\ref{thm main}
the
auxiliary
process is a percolation model 
similar to the Poisson blob model.
We call this auxiliary model  
totally asymmetric discrete Boolean
percolation. 
The comparison is not carried out via renormalization,
but rather the 
connected components in the auxiliary
percolation process
represent regions of space
traversed quickly,
while the vacant regions are traversed slowly.
The totally asymmetric discrete Boolean
percolation model and its properties
are described in Section \ref{sec TADBP percolation}.
In the proof of Theorem 
\ref{thm linear growth d geq 2}
the auxiliary process
is the greedy lattice animals model 
\cite{Martin02, GK94, CGGK93}.
Here too low values (in particular, zero)
in the greedy lattice animal model
represent sites
that do not have any quick outgoing
particles.

The idea to
get some information
about the spread of the process
by treating certain regions of space as
fast
appears in \cite{GM08},
where it is applied to 
a continuous-time continuous-space model of growing sets
introduced by Deijfen \cite{Dei03}.
To deal with the continuous-space nature of
Deijfen's model,
the authors in \cite{GM08}
introduce a continuous greedy paths model
which is 
a continuous-space
equivalent of the greedy lattice animals.
In the present paper we treat a lattice model,
hence we work directly with the greedy lattice animals.
Further discussion 
of the ideas of the proof 
can be found in 
Section \ref{sec proofs ideas}.

\begin{rmk}\label{liaise}
	Series 
	\eqref{obtuse}
	and \eqref{enclave}
	have independent convergent properties (that is,
	the convergence or divergence of either one of them
	does not imply anything about  the other).
	See
	Proposition \ref{independent series}
	for more details and examples.
\end{rmk}

\begin{rmk}
	We see in Theorem \ref{thm linear growth d = 1}
	and Theorem \ref{thm linear growth d geq 2}
	that it is possible for $\mu$
	to have infinite expectation 
	while the speed 
	is finite. 
	This may be considered counter-intuitive.
	One heuristic explanation for this may be
	that the probabilities
	of a simple random walk traveling 
	at high speed decline exponentially with 
	the distance (see Lemma \ref{guff}), and hence the conditions
	for the linear spread
	are given 
	in terms of 	roughly speaking the logarithmic moments.	 
\end{rmk}

Throughout the paper we use the following notation.
Let  \label{Page label notation}
$\{S_t, t \geq 0 \}$ be a simple continuous-time random walk
on $\Z ^\d$
and  $\tau _1, \tau _2, \dots$
be its jump times, $\tau _0 = 0$.
Let also $\{ (S _{t} ^{(x ,j)}, t \geq 0 ), x \in \Z ^\d, j \in \N \}$
be
independent copies of
$\{S_t, t \geq 0 \}$
assigned to individual particles.
For fixed $t, x$, and $ j$,
$x + S _{t} ^{(x ,j)} $ represents the position 
of  {$j$-th} particle started at location $x$,
$t$ units of time after the particle was activated.
For each realization of $\eta$, only the walks
$(S _{t} ^{(x ,j)}, t \geq 0 )$
with the indices satisfying $j  \leq  \eta (x)$
are used.
A
particle is identified with its index $(x, j)$.

\begin{rmk}
	In this work we always start with active particles 
	located exclusively 
	at the origin. Starting from
	a finite collection of sites with active particles
	does not affect the asymptotic spread rate
	and our results  still apply.
	This is a consequence of the following observation.
	Let $\zeta$ be a finite non-empty subset of $\Z ^\d$
	and  denote by $\A ^{\zeta}_t$
	the set of sites visited by the time $t$
	if at time $0$
	the locations of active particles 
	are exactly  $\zeta$.
	In the case $\eta(x) = 0$
	for some $x \in \zeta$, 
	an active particle is added to $x$ at time $0$ 
	(the addition of new particles is not necessary
	if active particles exist elsewhere.
	It is done for convenience
	because with the addition we get equality in \eqref{denuded}; 
	otherwise we would have to work with inclusions).
	Then 
	for any finite $\zeta \subset \Z ^\d$
	a.s.
	\begin{equation}\label{denuded}
	\A ^{\zeta}_t = \bigcup\limits _{x \in \zeta} \A ^{\{x\}}_t,
	\end{equation}	
	and all the conclusions about the linear or superlinear 
	spread rate follow.
	
\end{rmk}

\begin{rmk} \label{rmk frog1}
	It was shown in
	\cite{frog1}
	that
	the set $\A _t$
	can become infinite in a finite time
	if the tails of $\mu$ are heavy enough.
	In the present paper we address
	the conditions for the
	linear and superlinear spread rates.
	The observation in  \cite{frog1}
	raises the questions
	about the conditions
	separating the case of the superlinear spread rate
	such that at every moment of time $t > 0$
	only finitely many sites have been visited by active particles,
	and the case of 
	an explosion. 
	By the explosion here we mean
	that by a certain finite time,
	infinitely many sites have been  visited
	by active particles. 
	Further discussion can be found in
	a recent preprint \cite{frogE}.
\end{rmk}

\begin{rmk}
	As mentioned in Remark \ref{rmk frog1},
	the set $\A _t$
	can become infinite in a finite time
	when the tails of $\mu$ are heavy enough.
	It  therefore  behooves us
	to say a few words about the construction of the process. 
	Define the explosion time
	\begin{equation}
	\tau _e = \sup\{t: \A _t \text{ is finite}  \}.
	\end{equation}
	Prior to $\tau _e$ only finitely many events occur,
	hence the construction 
	on $[0,\tau _e)$ presents no challenges
	(indeed, on $[0,\tau _e)$
	the collection of active particles
	can be seen as a pure jump type Markov process, see e.g. \cite[Section 12]{KallenbergFound}).
	On the event $\{\tau _e < \infty \}$
	the construction of the process on $[\tau _e, \infty )$
	may present additional challenges
	because the process
	might become a system of  infinitely many interacting particles.
	Since we are only interested in the spread rate, 
	there is no need to consider the process on $[\tau _e, \infty )$;
	since $\A _t$ is non-decreasing in $t$
	and $\bigcup_{t < \tau _e} \A _t$ is infinite,
	we define
	the spread to be superlinear 
	on the event $\{\tau _e < \infty \}$.

\end{rmk}

\begin{rmk}\label{rift}
	Items (i) and (iii) of Theorem \ref{thm main}
	are statements of the form `if for some $B>1$,
	the series \dots converges/diverges, then \dots'
	For series 
	in \eqref{wack} and \eqref{enclave},
	the convergence for every $B>1$ is equivalent to
	the convergence for some $B > 1$. For instance 
	for
	if $1 < A < B$, then 
	\begin{equation*}
	\sum\limits _{m \in \N} \left[
	\mu(  [ B^{m} , \infty ) )
	\right] ^{\frac 1 \d}  \leq \sum\limits _{m \in \N} \left[
	\mu(  [ A^{m} , \infty ) )
	\right] ^{\frac 1 \d}
	\end{equation*}
	and 
	\begin{equation*}
	\sum\limits _{m \in \N} \left[
	\mu(  [ A^{m} , \infty ) )
	\right] ^{\frac 1 \d}
	\leq 
	\lceil \log _A B \rceil 
	+
	\lceil \log _A B \rceil 
	\sum\limits _{m \in \N} \left[
	\mu(  [ B^{m} , \infty ) )
	\right] ^{\frac 1 \d} 
	\end{equation*}
	Similarly for \eqref{enclave}
	\begin{equation*}
	\sum\limits _{n \in \N}   \mu \left( [ A ^ {n \ln ^2 n}, \infty) \right) 
	\simeq 
	\sum\limits _{n \in \N}   \mu \left( [ B ^ {n \ln ^2 n}, \infty) \right)
	\end{equation*}
	($\simeq $ means `have the same convergence properties'
	and is introduced in Section \ref{sec Notation and conv}).
	This is however not the case with the series in \eqref{obtuse}
	which may have different convergent properties
	for different $B>1$.
	
\end{rmk}

\subsection{Totally asymmetric discrete Boolean percolation (TADBP)}
We now describe the auxiliary percolation model
used in the proofs of Theorem \ref{thm linear growth d = 1}
and 
(ii) and (iii)
of
Theorem
\ref{thm main}.
It belongs to the class of discrete Boolean percolation.
An overview of the earlier works 
related to this class of models
can be found in 
\cite{BMS05},
while some connectivity properties
are established in \cite{CMG20}.
Here we are interested in the
case when the random connected neighborhoods,
or grains in the terminology of \cite{BMS05},
are totally asymmetric in the sense
that instead of the random balls $[x-r, x+r]$
with  random radii $r$,
the
intervals $[x,x+r]$ 
comprise connected components.

Let $\{\psi _z\}_{z \in \Z}$ be a collection of
independent identically distributed
$\Z _+$-valued random  variables
with distribution  $p_k = \PP{\psi _0 = k}$. 
We say that $x , y \in \Z$, $x \leq y$, are directly connected  (denoted by $x \xrightharpoondown{\Z}  y$)
if there exists $z \leq x$, $z \in \Z$, such that $z + \psi _z \geq y$.
We say that $x$  and $y$ are connected (denoted by $x \xrightarrow{\Z} y$)
if they are directly connected, 
or
if there exists $z_1 \leq ... \leq z_n \in \Z$,
$z_1 \leq x$, $z_n \leq y$, such that $x \in [z_1, z_1 + \psi _{z_1}]$, 
$y \in [z_n, z_n + \psi _{z_n}]$,
and $z_{j+1} \in [z_j, z_j + \psi _{z_j}]$ for $j = 1,2, ..., n-1$.
For a subset $Q \subset \Z$,
$x \xrightharpoondown{Q}  y$
and 
$x \xrightarrow{Q} y$
defined in the same way with an additional requirement 
that $x, y, z,  z_1, ..., z_n \in Q$ (in this paper we only consider 
$Q = \Z$ and $Q = \Z _+$).
The set $\Z$ is split into connected components.
We say that $x \in \Z$
is \emph{wet} if
the interval $[x-1,x]$
is contained in $[y, y + \psi _y]$
for some $y \in \Z$.
In other words, 
$x \in \Z$
is {wet} if for some $y\in \Z $,  $y < x$, $y + \psi _y \geq x$.
The sites that are  not wet are said to be dry. 
Note that
the statement
`$x$ is isolated'
is equivalent to 
`$x$ is dry and $\psi _x = 0$',
and that $x$ is wet if and only if $x-1$ and $x$ are connected.
We call the resulting random structure
totally asymmetric discrete Boolean
percolation (TADBP).
When considering TADBP on $\Z _+$,
we also talk about `wet' sites, 
with the understanding that both $x$ and $y$
are required to be from $\Z_+$.
Also,  
we consider the origin
to be wet for  TADBP on $\Z_+$.

The TADBP model on $\Z_+$
probably appears first
in \cite[Section 3]{Lamp70}, 
where the following interpretation is given.
At each $x \in \Z_+$ there is  a fountain
capable of wetting the sites ${x+1, \dots, x + \psi _x}$.
The fountain does not wet its own site.
In the case $\psi _x = 0$
the fountain at $x$ wets no site; it fails to operate. 
The model is very similar to the Poisson Boolean percolation model, or Poisson blob process, (\cite{MR96, FM07}), with the main differences being
the asymmetric nature of the random  sets
around each point
and
that the points of the Poisson point process in the Poisson Boolean  model
are replaced with the set of integers.
The model is further discussed in Section \ref{sec TADBP percolation}. 
A continuous-space equivalent of TADBP
is treated in \cite{TABP},
where the results analogous to those
in Section \ref{sec TADBP percolation}
are given.

Let $Q = \Z_+$ or $Q = \Z$.
Define the events $ x \xrightarrow{Q} \infty $ 
as that
for every $n \in \N$,
$x \in Q$
is connected to  $x + n$: $x \xrightarrow{Q} x + n$.
Note that if $p_0 = 0$, then trivially
$\P \{ x \xrightarrow{\Z} \infty  \} = \P \{ y \xrightarrow{\Z_+} \infty  \} = 1$
for every $x \in \Z$, $y \in \Z_+$.

\subsection{Very brief outlines of the proofs}\label{sec proofs ideas}

Item (i) of
Theorem \ref{thm main}
is split up into Theorems 
\ref{thm linear growth d = 1}
and \ref{thm linear growth d geq 2}.
Here we discuss the ideas of the proof
of (ii) and (iii)
of
Theorem \ref{thm main},
Theorem
\ref{thm linear growth d = 1},
and
Theorem
\ref{thm linear growth d geq 2}.	
Recall that $\{ (S _{t} ^{(x ,j)}, t \geq 0 ), x \in \Z ^\d, j \in \N \}$ 
are the random walks assigned to individual particles, and were introduced
on Page \pageref{Page label notation}. 

\emph{The ideas of the proofs of Theorem \ref{thm linear growth d = 1} and 
	(ii)
	of
	Theorem
	\ref{thm main}}.
The ideas  discussed here
are inspired by those articulated in \cite{GM08}
and applied there to Deijfen's model. 
Let $\d = 1$.
For $x \in \Z$
and $A > 0$
define 
\begin{equation}\label{secrete early}
\ell ^{(A)} _x = \max \Big\{ k \in \Z _+:
\exists t >0, j \in \overline{1, \eta (x)} \text{ such that } 
\frac{S _{t} ^{(x ,j)} }{t} \geq A \text{ and } S _{t} ^{(x ,j)} \geq k  \Big\} { \vee 0}.
\end{equation}
If for a given $x$ and $A$ the set on the right hand side of \eqref{secrete early}
is empty, then $\ell ^{(A)} _x = 0$. 
The random variable
$\ell ^{(A)} _x$
can be thought of as the length
of the longest interval
traveled toward $+\infty$
at an average speed at least $A$ 
by a particle started from $x$.

Consider now   TADBP with $\psi _x = \ell ^{(A)} _x$. 
The dry sites can be thought of
as being traveled over at a low speed below $A$,
whereas the wet sites are traveled over at a
high speed, at least $A$.
Imagine that $x \xrightarrow{\Z_+} \infty$. 
It means that for every $y \in (x, \infty)$
there exists $z< y$ such
that there is a particle starting 
from $z$ and traveling to $y$
or farther at speed at least $A$. 
Thus, intuitively, the speed
of the system should be at least $A$. If this is true 
for any $A>1$, then the spread must be superlinear.
Conversely, imagine many
sites of $\Z _+$
are dry. Then each of those sites
is traveled at speed not greater than $A$.
If such sites constitute a positive proportion of all sites
in a certain sense,
then we get a bound on the speed, and thus
the spread is linear. 

\emph{The ideas of the proof of  (iii)
	of
	Theorem
	\ref{thm main}}.
If \eqref{obtuse} converges, 
then we are under the assumptions in 
(ii)
of
Theorem
\ref{thm main}.
If it diverges, then 
the TADBP
random variables $\ell ^{(A)} _x$
may be
in the setting of Lemma \ref{profligacy}.
That is, the average size of a connected component
is infinity, but no site is connected to $+\infty$.
Thus, almost all sites are wet,
but the set of  dry sites
is unbounded. 
To deal with the dry sites,
we use a bound on the differences between activation time 
in some ways 
similar to Lemma 5.2 in \cite{stocCombust}. 
Let $\sigma _x = \min\{t \geq 0: x \in \A _t  \}$
be the moment when $x$ is visited by an active particle
for the first time.
We show that 
for sufficiently large $q$
\begin{equation}
\PP{ \sigma_{x} - \sigma _{x - 1} \geq q }
\leq  c_\varepsilon ^{q^{1/2}} .
\end{equation}
where $c_ \varepsilon     < 1$,
and then proceed to obtain
\begin{equation}\label{morose early}
\P \Big\{   \max\limits_{1 \leq y \leq x} (\sigma_{y} - \sigma _{y - 1}) \geq  C^2 {\ln ^2 x} \text{ infinitely often} \Big\} = 0
\end{equation}
for some constant $C > 0$. 
Combination of \eqref{enclave}
and
\eqref{morose early}
is then shown to imply the superlinear spread. 
The dry sites are `slow' and are dealt with using \eqref{morose early}; the wet sites are traveled at 
a speed at least $A$.

\emph{The ideas of the proof of Theorem \ref{thm linear growth d geq 2}}.
Here the greedy lattice animals \cite{Martin02, CGGK93}
play the role of the auxiliary model
instead of  TADBP. 
Recall that 
$\ell ^{(A)} _x$ were defined in 
\eqref{secrete early}, 
where the interpretation of 
$\ell ^{(A)} _x$ is also briefly discussed. 
Now, imagine that 
for any infinite sequence  $x_0 = \0, x_1, x_2, \dots, x_n, \dots$, $x_i \in \Z ^\d$, $x _i \ne x_j$, $i \ne j$,
$\min\limits _{j = 1,...,i-1}|x_{i} - x_j| = 1$, $i = 1,2,...,n$,  
with distinct points
the inequality 
\begin{equation*}
\frac{1}{n}\sum\limits _{i = 0} ^{n}
\ell _{_A} ^{(x_i)} \leq \frac 12
\end{equation*}
holds.
That means along the path $x_0 = \0, x_1, x_2, \dots $
not more than half the distance is traveled at speed greater than $A$.
The remaining one half is then traveled at speed at most $A$.
If this is true uniformly across
all paths, the linear spread should follow. 

\begin{rmk}
	The conditions imposed 
	on the sequence  $x_0 = \0, x_1, x_2, \dots, x_n, \dots$
	ensure that for every $m \in \N$
	the set $\{ x_0, x_1,...,x_m\}$
	is connected, however it is not necessary that $x_i$
	and $x_{i+1}$, $i \in \N$, are neighbors. This choice is 
	due to
	our desire
	to couple the stochastic combustion growth
	process with  the greedy lattice animals,
	and the way the latter are defined
	\cite{Martin02, CGGK93}.
\end{rmk}

\subsection{Notation and conventions}\label{sec Notation and conv}

For two series $\sum\limits _n a_n $ and $\sum\limits _n b_n $
with non-negative elements
we write $\sum\limits _n a_n \simeq  \sum\limits _n b_n  $
\label{Page simeq}
if they have the same convergence properties, that is,
they either both converge or both diverge. Respectively, 
we write $\sum\limits _n a_n \precsim  \sum\limits _n b_n  $
if  $\sum\limits _n b_n $ diverges, or if 
both  $\sum\limits _n a_n $ and
$\sum\limits _n b_n $ converge.
This is true for example if $a_n \leq b_n$
for large $n \in \N$ (but not necessarily for all $n \in \N$).

The minimum  and maximum operators $\wedge$ and $\vee$
precede addition and subtraction but follow after multiplication and division;
in other words, $a + b \vee cd =a + (b \vee (cd)) $.
For an interval $\I$, $|\I|$ is its length.
We adopt the following convention 
regarding the operations over 
the empty set:
$\sum\limits _{q \in \varnothing}  q = 0$, 
$\prod\limits _{q \in \varnothing}  q = 1$,  $\bigcup\limits _{q \in \varnothing}  q = \varnothing$, 
$\max \varnothing = \sup \varnothing = - \infty$, 
$\min \varnothing = \inf \varnothing = + \infty$.
The symbol $\1$ denotes an indicator.

\begin{rmk}
	The following inequalities are used extensively in the paper. 
	For $a \in(0,1), b \geq 1$,
	\begin{alignat*}{3}
	(1 - a) ^b & \geq 1 - 1 \wedge ab, \ \ \  \ \ \ \ \ \ 
	(1 - a) ^b & & \leq e^{-ab}
	\\
	1 - (1 - a) ^b & \geq 1 - e^{-ab}, \ \ \  \ \
	1 - (1 - a) ^b & & \leq 1 \wedge ab.
	\end{alignat*}
	They are consequences and extensions of Bernoulli's inequality.
\end{rmk}

\section{Totally asymmetric discrete Boolean percolation: properties}\label{sec TADBP percolation}

In this section we
determine the fraction 
of wet sites and
establish
necessary and sufficient conditions
for 
a node being connected
to $+\infty$
with positive probability.
Most of the
results in this section
are not new. In particular,
the transience criterion
in
Proposition \eqref{saunter}
appears in  Kesten's Appendix to \cite{Lamp70}
and also later in \cite{Kel06}, and the positive recurrence criterion
is formulated in \cite[Page 283]{Kel06};
it is also the content of \cite[Proposition 1.1]{Zer18}.
The proof of our Lemma \ref{profligacy} is basically the same as the proof
of Proposition 1.1 in \cite{Zer18}.

First we establish under what conditions $\P \{ x \xrightarrow{\Z} \infty  \} = 1$.
For $ n \in \N$, let $r_n = \sum\limits _{i = n} ^\infty p_i = \PP{ \psi _x \geq n}$
be the tail of the distribution $\{p _n \}_{n \in \Z _+}$.
Let the exclamation mark 
in front of the connectivity relations
denote the negation,
for example 
the event
$ \{ -1 \ ! \! \!  \xrightarrow{\Z} 0  \}$
is the complement of $\{ -1 \xrightarrow{\Z} 0  \}$.
We exclude trivial cases and assume $p_0 \in (0,1)$.

\begin{lem}\label{compel}
	Consider  TADBP on $\Z$.
	We have $\P \{ x \xrightarrow{\Z} \infty  \} \in \{ 0,1\}$,
	and  $\P \{ x \xrightarrow{\Z} \infty  \} = 1$
	if and only if 
	$\prod\limits_{  k =1 } ^\infty (1 - r_k) = 0$.
	Respectively, 
	$\P \{ x \xrightarrow{\Z} \infty  \} = 0$
	if and only if 
	$\prod\limits_{  k =1 } ^\infty (1 - r_k) > 0$,
	and in this case a.s.
	\begin{equation}\label{countenance}
	\frac{\# \{ k \in \{  1, \dots, n \}: k  \text{\emph{ is dry}} \}}{n} \to 
	\prod\limits_{  k =1 } ^\infty (1 - r_k), 
	\ \ \ n \to \infty.
	\end{equation}
\end{lem}
\begin{proof}
	We have 
	\begin{equation}\label{pug}
	\P \{ -1 \ ! \! \!  \xrightarrow{\Z} 0  \} =
	\PP{ \psi _{-m} < m \text{ for all } m \in \N  }
	=\prod\limits _{m = 1} ^\infty
	\PP{ \psi _0 < m   }
	=
	\prod\limits_{  m =1 } ^\infty (1 - r_m).
	\end{equation}
	Hence $\P \{ -1  \xrightarrow{\Z} 0  \} = 1$ if $\prod\limits_{  k =1 } ^\infty (1 - r_k) = 0$,
	and by translation invariance  $\P \{ x  \xrightarrow{\Z} x+1  \} = 1$
	for every $x \in \Z$. Thus, a.s. every node is connected
	to infinity provided $\prod\limits_{  k =1 } ^\infty (1 - r_k) = 0$.
	
	Let now $\prod\limits_{  k =1 } ^\infty (1 - r_k) > 0$. 
	Define the random variables 
	$Z_n = \1 \{ n-1   ! \! \!  \xrightarrow{\Z} n  \}$.
	Since $\{\psi _n\}_{n \in \Z}$
	is a sequence of i.i.d. random variables and thus ergodic,
	so is $\{Z _n\}_{n \in \Z}$, because $$Z_n = \prod\limits_{  k =1 }\1\{\psi_{n-k} < k  \} $$
	is a functional transformation of $\{\psi _n \}$ (see 
	e.g. \cite[Theorem 7.1.3]{Durr_Theory_and_examples}).
	By the ergodic theorem
	and \eqref{pug}
	a.s.
	\begin{equation}
	\frac{\sum\limits _{k = 1} ^n Z_k}{n} \to \E Z_1 = \prod\limits_{  k =1 } ^\infty (1 - r_k) > 0.
	\end{equation}
\end{proof}

In the remaining part of the section
we focus on  TADBP on $\Z_+$.

\begin{lem}\label{quail}
	Consider  TADBP  
	on $\Z _+$ and let $\prod\limits_{  k =1 } ^\infty (1 - r_k) >0$.	
	The fraction of sites that are dry is 
	$\prod\limits_{  k =1 } ^\infty (1 - r_k) $
	in the sense that a.s.
	\begin{equation} \label{slovenly}
	\frac{\# \{ k \in \{  1, \dots, n \}: k  \text{\emph{ is dry }} \}}{n} \to 
	\prod\limits_{  k =1 } ^\infty (1 - r_k), 
	\ \ \ n \to \infty.
	\end{equation}
	The fraction of isolated sites is 
	$p_0 \prod\limits_{  k =1 } ^\infty (1 - r_k) $.
	A.s.
	no site is connected to $+\infty$.
\end{lem}
\begin{proof}
	The convergence \eqref{slovenly} follows from \eqref{countenance}
	since a.s. $\sup\limits _{m \in \N} (\psi _{-m} - m) < \infty$,
	and for sites $x > \sup\limits _{m \in \N} (\psi _{-m} - m) $
	the values $\psi _{-m}$, $m \in \N$, do not
	have any effect on whether $x$
	is wet or dry. 
	
	A site $x$ is isolated if and only if $x$ is  dry
	and $\psi _x = 0$. Since $\psi _x$
	is independent of $\{\psi _y \}_{y < x}$, the second statement of the lemma follows. 
\end{proof}

For $m \in \Z _+$, denote by $Y_m$ the difference between 
the rightmost site directly connected to $m$, and $m$, that is,
\begin{equation}\label{commandment}
Y_m = \max\{ l : m \xrightharpoondown{\Z _+} l  \}   - m.
\end{equation}
Recall that by definition  $ m\xrightharpoondown{\Z} m$
holds true for all $m \in \Z _+$, hence $Y_m \geq 0$.
By construction for $m \in \N$
\begin{equation}
Y_m =  \psi _m \vee (Y_{m-1} - 1)  = \psi _m \vee (\psi _{m-1} - 1)
\vee ... \vee (\psi _1 - m + 1) \vee (\psi _ 0 - m).
\end{equation}
Note that since we have assumed $p_0 \in (0,1)$,
$(Y_t, t  \in \Z_+)$
constitutes an irreducible Markov chain on $\Z _+$.
In essence
this Markov chain appears in \cite[(24)]{Lamp70}
and \cite[Page 268]{Kel06}.

\begin{lem} \label{goad}
	Assume that 
	\begin{equation} \label{onerous}
	\sum\limits _{n = 1} ^ \infty \prod_{k = 1}^{n} (1 - r_k)< \infty.
	\end{equation}	 
	Then a.s. there exists  $x \in \Z _+$
	connected to $\infty$:
	\begin{equation} \label{sycophancy}
	\PP{  x \xrightarrow{\Z _+} \infty  \text{ for some } x \in \N } = 1.
	\end{equation}

\end{lem}
\begin{proof}
	We have
	\begin{equation} \label{snazzy}
	\PP{ Y _m  = 0} = \prod_{i = 0}^{m} \PP{ \psi _i \leq m - i } = \prod_{i = 0}^{m} \PP{ \psi _i \leq i }
	=
	\prod_{i = 0}^{m} (1 - r _{i + 1}).
	\end{equation}
	By \eqref{onerous} and \eqref{snazzy}
	\begin{equation}
	\PP{ Y _m  = 0 \text{ for infinitely many } m \in \N} = 0.
	\end{equation} 
	It remains to note that if
	some $m \in \N$,
	$Y_i \geq 1$ for $i \geq m$,
	then $m \xrightarrow{\Z _+} \infty$.
\end{proof}

\begin{lem} \label{profligacy}
	Assume that 
	\begin{equation} \label{grovel}
	\sum\limits _{n = 1} ^ \infty \prod_{k = 1}^{n} (1 - r_k) =  \infty
	\end{equation}	 
	and 
	\begin{equation} \label{jaunty}
	\prod_{k = 1}^{\infty} (1 - r_k) =  0.
	\end{equation}	
	Then 
	\begin{equation} \label{fight tooth and nail}
	\PP{  x \xrightarrow{\Z _+} \infty  \text{ for some } x \in \N } = 0.
	\end{equation}

\end{lem}
\begin{proof}
	Recall that the Markov $(Y_t, t  \in \Z_+)$
	is defined in \eqref{commandment}.
	Since 
	\begin{equation}
	\{  x \xrightarrow{\Z _+} \infty   \} = \{ Y_k >0, k = x, x+1, x+2, \dots \},
	\end{equation}
	\eqref{fight tooth and nail} is equivalent 
	to $(Y_t, t  \in \Z_+)$ being recurrent.
	By \eqref{grovel}
	\begin{equation}
	\sum\limits _{n = 1}\P\{Y_n = 0 \} = 
	\sum\limits _{n = 1} ^ \infty \prod_{k = 1}^{n} (1 - r_k) =  \infty.
	\end{equation}
	Therefore $(Y_t, t  \in \Z_+)$ is recurrent
	by the well known properties of Markov chains
	with a countable state space,
	see e.g. \cite[Theorem 6.4.2]{Durr_Theory_and_examples}
	or \cite[Theorem 1, Section 5, Chapter 8]{ShiryaevProb-2}.
\end{proof}

\begin{rmk}
	Lemma \ref{profligacy} is akin to the dichotomy
	occurring under certain conditions
	in Boolean percolation
	when each occupied component 
	is a.s. finite
	but the expected size
	of an occupied component is infinite, 
	see  \cite[Corollary 3.2]{MR96} 
	or \cite{TABP}.
\end{rmk}

We note that the individual assumptions 
of
Lemmas \ref{quail}, \ref{goad} and   \ref{profligacy}
regarding $\{p_i\}_{i \in \Z_+}$,
\begin{equation}\label{high jinks}
\prod\limits_{  k =1 } ^\infty (1 - r_k) >0,
\end{equation}
\begin{equation}\label{high jinks2}
\sum\limits _{n = 1} ^ \infty \prod_{k = 1}^{n} (1 - r_k) < \infty,
\end{equation}
and
\begin{equation}\label{high jinks3}
\sum\limits _{n = 1} ^ \infty \prod_{k = 1}^{n} (1 - r_k) =  \infty,
\  \  \ 
\prod_{k = 1}^{\infty} (1 - r_k) =  0.
\end{equation}	
exhaust all options; that is, one (and only one of course) of the conditions 
\eqref{high jinks}, \eqref{high jinks2},
and \eqref{high jinks3} always holds.
Therefore, those lemmas
lead to characterization of the recurrence 
properties of the Markov chain $(Y_t, t  \in \Z_+)$.
These properties are collected in the next proposition.
It is formulated in the self-sufficient way,
so that all necessary notation used in this
section is reintroduced. As indicated
at the beginning of the section, 
this proposition does not contain new results. 

\begin{prop}\label{saunter}
	Let $\{ \psi _k\}_{k \in \Z_+}$ be 
	a sequence of $\Z_+$-valued
	random variables with distribution $\{p_i\}_{i \in \Z_+}$.
	Set  $r_k = \sum\limits _{i = k} ^\infty p_i$
	and 
	$Y_m =  \psi _m \vee (Y_{m-1} - 1)$, $m \in \N $, $Y_0 = \psi _0$.
	Then 
	\begin{itemize}
		\item [] (i) The Markov chain $(Y_t, t  \in \Z_+)$
		is positive recurrent if and only if $\ \E \psi _1 = \sum_{k = 1}^{\infty} r_k  < \infty$,
		\item [] (ii) $(Y_t, t  \in \Z_+)$
		is transient if and only if  
		$ \
		\sum\limits _{n = 1} ^ \infty \prod\limits_{k = 1}^{n} (1 - r_k)< \infty,
		$	
		\item [] (iii)  The chain $(Y_t, t  \in \Z_+)$  is null recurrent
		if and only if
		both
		$\ \sum_{k = 1}^{\infty} r_k  = \infty$
		and 	$ {
			\sum\limits _{n = 1} ^ \infty \prod\limits_{k = 1}^{n} (1 - r_k)=\infty }
		$.
		
	\end{itemize}
\end{prop}

\begin{proof}
	As was noted in the proof of Lemma \ref{profligacy},
	\eqref{fight tooth and nail} is equivalent 
	to $(Y_t, t  \in \Z_+)$ being recurrent. 
	Thus
	Lemmas  
	\ref{quail}, 
	\ref{goad},
	and \ref{profligacy}
	combined yield (ii).
	
	Assume now 	$ \
	\sum\limits _{n = 1} ^ \infty \prod_{k = 1}^{n} (1 - r_k)=\infty
	$,
	that is, that the chain is recurrent. 
	If $ \sum_{k = 1}^{\infty} r_k  = \infty$,
	then the chain is cannot be positive recurrent because the 
	expected recurrence time to $0$ is
	greater than $\sum_{k = 1}^{\infty} r_k$, so it is infinite. 
	On the other hand,
	if $\sum_{k = 1}^{\infty} r_k < \infty$,
	then for every $m \in \N$
	\begin{multline*}
	\P \{Y _m = 0\} = \PP{ \psi _m \leq 0} \PP{ \psi _{m-1} \leq  1}
	\cdot \ldots \cdot \PP{ \psi _{0} \leq  m} 
	=  \prod_{k = 1}^{m+1} (1 - r_k) 
	\geq \prod_{k = 1}^{\infty} (1 - r_k)  > 0.
	\end{multline*}
	and hence $(Y_t, t  \in \Z_+)$
	cannot be null recurrent.
\end{proof}

In the null-recurrent case a.s. 
$$\frac{\#\{ k \in \{1,...,n\}: Y_k  =0 \}}{n} \to 0, \ \ \ n \to \infty,$$
and since $k$ is dry if and only if $Y_k = 0$, we get   
\begin{cor}
	Assume that \eqref{grovel} and \eqref{jaunty}  hold. Then
	for TADBP on $\Z_+$
	the fraction of wet 
	sites is $1$.
\end{cor}

The next two lemmas will be helpful 
in showing that the spread is superlinear.
They provide a way to translate
the properties
of the associated
TADBP to the properties
of the stochastic combustion growth process. 

\begin{lem}\label{hold feet to the fire}
	Let $x \in \N$. A.s. on $\{ x \xrightarrow{ \Z _+ } \infty \}$,
	every site $y > x$ is wet, and there exists
	a (random) sequence $x = x_0 < x_1 < x_2 < \dots $, $x_i \in \N$,
	such that for every $i \in \Z _+$
	\begin{equation}\label{slink}
	x _{i + 1} \leq x _ i + \psi _{x _ i} <  x _{i + 2}.
	\end{equation}
	In particular, every $z \geq x$  belongs to no more than two
	intervals of the type $[x_i, x_i + \psi _{x _ i} ]$, $i \in \Z _+$.
\end{lem}
\begin{proof}
	By definition of $\xrightarrow{ \Z _+ }$, 
	every site $y > x$ is wet a.s. on $\{ x \xrightarrow{ \Z _+ } \infty \}$.
	Define the elements of the sequence $\{ x_i\}_{i \in \Z_+}$
	consecutively setting $x_0 = x$
	and letting for $i \in \Z_+$
	\begin{equation}
	x_{i+1} = \max\{y \in [x_i + 1, x_i + \psi _{x _ i} ] \cap \N : y + \psi _y = \max\{z + \psi _z: z = x_i + 1, \dots, x _i + \psi _{x _i} \}   \}.
	\end{equation}
	In other words, $x_{i+1} \in [x_i + 1, x_i + \psi _{x _ i} ]$ is 
	characterized by two properties:
	\begin{itemize}
		\item[(i)]
		for every 
		$z \in [x_i + 1, x_i + \psi _{x _ i} ] \cap \N  $, 
		\begin{equation*}
		x_{i+1} + \psi _{x _{i+1}} \geq z + \psi _z,
		\end{equation*}
		\item[(ii)]
		and for every $z' \in [x_{i+1} + 1, x_i + \psi _{x _ i} ] \cap \N  $,
		\begin{equation*}
		x_{i+1} + \psi _{x _{i+1}} > z' + \psi _{z'}.
		\end{equation*}
	\end{itemize}
	(here $[a,b] = \varnothing$ if $a > b$).
	By construction, $x _{i + 1} \leq x _ i + \psi _{x _ i}$, so 
	the left inequality in \eqref{slink} holds.
	A.s. on $\{ x \xrightarrow{ \Z _+ } \infty \}$, 
	$x_{i+1} + \psi _{x _{i+1}} >  x_i + \psi _{x _ i}$, because otherwise
	$x_i + \psi _{x _ i} + 1$ would not be wet.
	Hence a.s. on $\{ x \xrightarrow{ \Z _+ } \infty \}$ also $x_{i+2} + \psi _{x _{i+2}} >  x_{i+1} + \psi _{x _{i+1}} $. 
	Therefore the inequality  $    x _{i + 2} \leq x _ i + \psi _{x _ i} $
	is impossible a.s.  on $\{ x \xrightarrow{ \Z _+ } \infty \}$ because it would contradict (i).
\end{proof}

The next lemma replicates Lemma \ref{hold feet to the fire}
for the case of a finite component. The proof is
practically identical
and is therefore omitted. 
\begin{lem}\label{hold feet to the fire2}
	Let $x \in \N$. A.s. on $\{ x \xrightarrow{ \Z _+ } y \}$,
	every site $z \in [x+1,y]$ is wet, and there exists
	a (random) sequence $x = x_0 < x_1  < \dots < x_m = y $, $x_i \in \N$,
	such that for every $i \in \{ 0, ..., m-2\}$
	\begin{equation}\label{slink2}
	x _{i + 1} \leq x _ i + \psi _{x _ i} <  x _{i + 2},
	\end{equation}
	and $	 x _ {m-1} + \psi _{x _ {m-1}} =  x _{m}$.
	In particular, every $z \in [x,y]$  belongs to no more than two
	intervals of the type $[x_i, x_i + \psi _{x _ i} ]$, $i \in \{ 0, ..., m-1\}$.
\end{lem}

\section{Proofs of Theorem \ref{thm linear growth d = 1} and 
	(ii) 
	of
	Theorem
	\ref{thm main}}  \label{sec proofs directly from TADBP}

Let $\tip{t}$ be the position of the rightmost active 
particle at time $t\geq 0$. Note that it is not necessarily 
true
that  $\tip{t} = \sup \A _t$ for  $t \geq 0$, 
because active particles can move back toward $-\infty$.
Let us introduce a total ordering $\prec$ on the set of indices
with $(x,i) \prec (y,j)$ if  $x < y$, or if $x = y$
and $i < j$.

Recall that $\sigma _x = \min\{t \geq 0: x \in \A _t  \}$
is 
the activation time for particles at
$x$.
At  time $t\geq 0$, let $X_0 = (x_0, i _0)$
be the particle with the smallest index
with respect to $\prec$
located at $\tip{t}$.
The particle $X_0$
depends on $t$, even if it is not reflected in the notation.
We note here that 
excluding the case $x_0 < 0$,
no particle located at $\tip{t}$
at time $t$
was activated before 
$X_0$.
Let $X_1 = (x_1, i _1)$ be the particle that activated $X_0$, 
and further 
define  recursively $X_{k + 1}$
as the particle that activated $X_k$,
until $X_\m = (0, i_ \m)$ for some $\m \in \N$.
Set $w _i = \sigma _{X_k}$, $k = 1,\dots,\m$.
Let us note right here 
that the sequence $\{ X_i, i = 1,\dots,\m\}$
depends on $t$.
Denote 
by $\W _k$,
$k = 1, \dots, \m$,
the interval 
$[ w_{\m - k + 1}, w_{\m -k }  ]$.
If $x_{\m - k + 1} <  x_{\m -k } $,
the expression
\begin{equation}
\frac{ (x_{\m -k } - x_{\m -k +1} )\vee 0}{ w_{\m -k } - w_{\m - k + 1}   }
=
\frac{ (x_{\m -k } - x_{\m -k +1} )\vee 0}{ |\W _k|  }.
\end{equation}
can be seen as the speed  at which the interval $[x_{\m - k + 1}, x_{\m -k } ]$ is traversed.
We  take non-negative part in the numerator
because the sequence $\{ x_i, i = 1,\dots,\m\}$ does not have to be non-increasing;
indeed, it is possible that active particles
from the origin travel leftward, activate a particle at $-k$, $k \in \N$,
and that particle then moves toward $+\infty$
very quickly overtaking every other active particle,
and becomes a leading particle for some time.

Recall that 
for $x \in \Z$
and $A > 0$, we have defined
\begin{equation}\label{secrete}
\ell ^{(A)} _x = \max \{ k \in \Z _+:
\exists t >0, j \in \overline{1, \eta (x)} \text{ such that } 
\frac{S _{t} ^{(x ,j)} }{t} \geq A \text{ and } S _{t} ^{(x ,j)} \geq k  \}.
\end{equation}
Let $\{p_ k ^{(A)} \}_{k \in \Z _+}$ be the distribution of $\ell ^{(A)} _x$,
$p_ k ^{(A)} = \PP{ \ell ^{(A)} _x = k }$. Let $r_ k ^{(A)}$
be the corresponding tail,
$r_k ^{(A)}= \sum\limits _{i = k} ^\infty p_i ^{(A)} = \PP{ \ell ^{(A)} _x \geq k}$.
Note that $p_ 0 ^{(A)} \in (0,1) $, 
since 
\begin{multline*}
1 - p_ 0 ^{(A)} = 
\PP{\ell ^{(A)} _x \geq 1} = \PP{\ell ^{(A)} _x \geq 1, \eta(x) \geq 1}
\\
\geq 
(1-\mu(0)) \P \Big\{\tau _1 < \frac{1}{A}, S_{\tau_1} = 1  \Big\}
=  \frac{ (1-\mu(0))}{2} \left(1 -  \exp\Big\{ -\frac 1A  \Big\}  \right) > 0
\end{multline*}
and for $n\in \N$ with $\mu(n) > 0$, 
\begin{equation*}
p_ 0 ^{(A)}
=
\PP{\ell ^{(A)} _x =0 } \geq  \PP{\ell ^{(A)} _x =0, \eta(x) = n}
\geq 
\mu(n) \left(  \P\Big\{\forall t > 0: \frac{S_t}{t} < A \Big\}  \right)^n
>0.
\end{equation*}

The proofs of the next two propositions
are relatively short,
however they
contain the key 
steps 
in comparing the stochastic combustion 
growth process with TADBP.
In
Proposition \ref{prop linear growth}
the linear growth is established,
therefore 
TADBP on $\Z$ is used in the proof because
we want our stochastic combustion growth process
to be dominated  by TADBP
in a certain sense
as we establish the linear spread.
In contrast,
TADBP on $\Z_+$
is used
in the proof of Proposition \ref{prop superlinear growth},
because 
we want the stochastic combustion growth process
to dominate TADBP,
and it turns out that 
even
considering only the particles
on $\Z_+$ is enough 
to ensure the superlinear spread.

\begin{prop} \label{prop linear growth}
	
	Suppose that for  some $A >0$
	\begin{equation}\label{charade}
	\prod\limits_{  k =1 } ^\infty (1 - r_k ^{(A)}) > 0.
	\end{equation}
	Then $\sup  \A _t$ grows linearly in time. 
\end{prop}
\begin{proof}
	Consider  TADBP on $\Z$
	with interval distribution
	$\{p_ k \}_{k \in \Z _+} = \{p_ k ^{(A)} \}_{k \in \Z _+}$.
	By Lemma \ref{compel}, the fraction of  sites that are dry  is 
	$\rmu =   \prod\limits_{  k =1 } ^\infty (1 - r_k) > 0$.
	Therefore with high probability
	at least $(1 - \varepsilon)\rmu \cdot \tip{t} $ sites among $1,2,\dots, \tip{t}$
	are traveled at a speed at most $A$, where $\varepsilon \in (0,10^{-1})$ is a small number. 
	Hence
	a.s. for large $t$,  $\sum\limits _{ k = 1} ^ \m |\W _k| \geq \frac{ (1 - \varepsilon) \rmu \cdot \tip{t}}{A}  $
	and
	\[
	\limsup _{ t \to \infty }
	\frac{\tip{t}}{t} = \limsup _{ t \to \infty } \frac{\tip{t}}{\sum\limits _{ k = 1} ^ \m |\W _k|} \leq 
	\frac{A}{(1 - \varepsilon) \rmu}.
	\]
	Since $\sup  \A _t { =} \sup\limits _{s \leq t} \tip{s}$, the statement of the proposition follows.
\end{proof}

\begin{prop} \label{prop superlinear growth}
	Suppose that for  all  $A >0$
	\begin{equation} \label{lapidation}
	\sum\limits_{  m =1 } ^\infty \prod\limits_{  n =1 } ^m (1 - r_n ^{(A)}) < \infty.
	\end{equation}
	Then 
	$$\liminf\limits _{t \to \infty} \frac{\tip{t}}{t} = \infty, $$
	and
	$\sup  \A _t$ grows faster than linearly in time.
	
\end{prop}
\begin{proof}
	Here  
	for $A >0$
	we
	consider  TADBP
	on $\Z_+$
	with interval distribution
	$\{p_ k \}_{k \in \Z _+} = \{p_ k ^{(A)} \}_{k \in \Z _+}$.
	By \eqref{lapidation} and Lemma \ref{goad}
	a.s.
	there is a large (random) number 
	$x_0 \in \N$ such that $x_0  \xrightarrow{\Z _+}  \infty$. 
	In particular, every site $x \in (x_0, \infty ) \cap \N$
	is wet. 
	By Lemma \ref{hold feet to the fire}, 
	there exists a sequence $x_0, x_1, x_2, \dots$
	such that $x_i < x_{i+1} \leq x_i + \ell ^{(A)} _{x_i}$, $i \in \Z _+$,
	and
	every site $y \in (x_0, \infty ) \cap \N$
	belongs to at most two intervals of the type $[x_i, x_i + \ell ^{(A)} _{x_i}]$, $i \in \Z _+$.
	Consequently for every $n \in \N$
	\begin{equation} \label{obdurate}
	\sum\limits _{i = 0} ^{n}  \ell ^{(A)} _{x_i}
	=
	\sum\limits _{i = 0} ^{n} \big| [x_i, x_i + \ell ^{(A)} _{x_i}] \big|  \leq 2(x_n + \ell ^{(A)} _{x_n} - x _0).
	\end{equation}
	Once $x_0$ is reached,
	let us  consider only particles that start at $x_i$
	and 
	travel to  $ x_i + \ell ^{(A)} _{x_i}$ at speed at least $A$.
	For $n \in \N$
	and $y_n: = x_n + \ell ^{(A)} _{x_n}$
	by \eqref{obdurate} 
	\begin{equation*}
	\sigma _{y_n} - \sigma _{x _0} \leq 
	\frac 1 A \sum\limits _{i = 0} ^{n-1} \big| [x_i, x_i + \ell ^{(A)} _{x_i}] \big|  
	\leq  \frac{2(y _{n} - x_0)}{ A }.
	\end{equation*} 
	Hence 
	\begin{equation}
	\limsup\limits _{t \to \infty} \frac{\tip{t}}{t}
	\geq
	\limsup\limits _{n \to \infty} \frac{ y_n}{\sigma _{y_n}}
	=
	\limsup\limits _{n \to \infty} \frac{x_0 +  (y_n - x _0)}{\sigma _{x _0} + (\sigma _{y_n} - \sigma _{x _0})}
	=
	\limsup\limits _{n \to \infty} \frac{ y_n - x _0}{ \sigma _{y_n} - \sigma _{x _0}}
	\geq \frac{A}{2}.
	\end{equation}
	Since $A>0$ can be arbitrary large, the statement of the proposition follows.
\end{proof}

In the remaining part of this section 
we establish 
that 
\eqref{charade}
is equivalent to \eqref{egregious},
whereas 
\eqref{lapidation}
is equivalent to \eqref{obtuse}.
The next three lemmas
dealing with convergence
properties of certain related series.  
Recall that the notation $\simeq$
for two series was introduced
on Page \pageref{Page simeq}.

\begin{lem}\label{putrid}
	For $A > 1$, 
	\begin{equation}
	\sum\limits_{  k =1 } ^\infty \mu (k)
	\sum\limits_{  n =1 } ^\infty  1 \wedge k A ^{-n} < \infty 
	\end{equation}
	if and only if 
	\begin{equation}
	\sum\limits _{n = 1}  \mu(n) \log n  < \infty.
	\end{equation}
\end{lem}
\begin{proof} 
	We have 
	\begin{multline*}
	\sum\limits_{  k =1 } ^\infty \mu (k)
	\sum\limits_{  n =1 } ^\infty  1 \wedge k A ^{-n} = 
	\sum\limits_{  k =1 } ^\infty \mu (k) 
	\sum\limits_{  1 \leq n \leq \log _A k  } 1
	+
	\sum\limits_{  k =1 } ^\infty \mu (k) 
	\sum\limits_{  n > \log _A k    } k A ^{-n} 
	\\
	\simeq \sum\limits_{  k =1 } ^\infty \mu (k)  \log _A k 
	+ 
	\sum\limits_{  k =1 } ^\infty \mu (k) 
	\frac{1}{1 - A^{-1}}
	= (\log  A) ^{-1} \sum\limits_{  k =1 } ^\infty \mu (k)  \log  k 
	+ 
	\frac{1}{1 - A^{-1}}.
	\end{multline*}
\end{proof}

\begin{lem}\label{anodyne}
	Suppose that for some $ A, B > 1 $ 
	\begin{equation}\label{sulky}
	1 - \sum\limits_{  k =0 } ^\infty \mu (k)  \left[ 1 - A  ^{-n} \right] ^k 
	\leq r _n 
	\leq 
	1 - \sum\limits_{  k =0 } ^\infty \mu (k)  \left[ 1 - B ^{-n} \right] ^k,
	\ \ \ n \in \N.
	\end{equation}
	Then   
	\begin{equation}\label{salacious}
	\prod\limits_{ n =1 } ^\infty (1 - r_n ) = 0
	\end{equation}
	if and only if 
	\begin{equation}\label{deleterious}
	\sum\limits _{n = 1}  \mu(n) \log n  = \infty.
	\end{equation}

\end{lem}
\begin{proof}
	Note that by \eqref{sulky}, $r_n < 1$, $n \in \N$.
	Since $(1 + a) ^ b \leq e^{ab}$ for $a \geq -1$, $b \geq 0$, 
	we have 
	\begin{equation} \label{gyrate}
	1 - \sum\limits_{  k =0 } ^\infty \mu (k)  \left[ 1 - A  ^{-n} \right] ^k 
	\geq   1 - \sum\limits_{  k =0 } ^\infty \mu (k) e^{-k A ^{-n}}
	= \sum\limits_{  k =0 } ^\infty \mu (k) \left[ 1 - e^{-k A ^{-n}} \right].
	\end{equation}
	By 
	\eqref{sulky}
	and
	\eqref{gyrate}
	and  since $ \inf\limits _{a > 0} \frac{1 - e ^{-a}}{1 \wedge a} > 0$
	and
	$ \sup\limits _{a > 0} \frac{1 - e ^{-a}}{1 \wedge a} =1$
	we get 
	\begin{multline} \label{keep tabs on}
	\sum\limits_{ n =1 } ^\infty  r_n \geq 
	\sum\limits_{  n =1 } ^\infty
	\sum\limits_{  k = 0 } ^\infty \mu (k) \left[ 1 - e^{-k A ^{-n}} \right]
	= \sum\limits_{  k = 0 } ^\infty \mu (k)
	\sum\limits_{  n =1 } ^\infty \left[ 1 - e^{-k A ^{-n}} \right]
	\simeq 
	\sum\limits_{  k = 0 } ^\infty \mu (k)
	\sum\limits_{  n =1 } ^\infty  1 \wedge k A ^{-n}.
	\end{multline}
	On the other hand, 
	since 
	by Bernoulli's inequality 
	$(1 - a) ^b  \geq 1 - 1 \wedge ab$, $a \in (0,1)$, $b \geq 1$, we have
	\begin{equation}
	1 - \sum\limits_{  k = 0 } ^\infty \mu (k)  \left[ 1 - B ^{-n} \right] ^k 
	\leq 
	1 - \sum\limits_{  k = 0 } ^\infty \mu (k)  \left[ 1 - 1\wedge k B ^{-n} \right]
	= \sum\limits_{  k = 0 } ^\infty \mu (k)  \left[ 1\wedge k B ^{-n}\right].
	\end{equation}
	Hence
	\begin{equation} \label{off the charts}
	\sum\limits_{ n =1 } ^\infty  r_n \leq 
	\sum\limits_{ n =1 } ^\infty
	\sum\limits_{  k =1 } ^\infty \mu (k)  \left[ 1\wedge k B ^{-n}\right]
	= \sum\limits_{  k =1 } ^\infty \mu (k) 
	\sum\limits_{ n =1 } ^\infty
	1\wedge k B ^{-n}.
	\end{equation}
	From \eqref{keep tabs on}, \eqref{off the charts}, and  Lemma \ref{putrid} 
	it follows that 
	\begin{equation} 
	\sum\limits_{ n =1 } ^\infty  r_n \simeq 
	\sum\limits_{  k =1 } ^\infty \mu (k)
	\sum\limits_{  n =1 } ^\infty  1 \wedge k A ^{-n}
	\simeq \sum\limits _{n = 1}  \mu(n) \log n.
	\end{equation}
	The equivalence of 
	\eqref{salacious} and \eqref{deleterious}
	follows since 
	\eqref{salacious} is equivalent to 
	$  \sum\limits_{ n =1 } ^\infty  r_n = \infty$.
\end{proof}

Regarding the series in \eqref{icky} in the following lemma, note that
\begin{equation}
\sum\limits _{m=1} ^\infty \prod\limits _{n = 1} ^ m 
\sum \limits _{k = 0} ^{\infty} \mu (k) \left[ 1 - B ^{-n}\right]^{k}
= \sum\limits _{m=1} ^\infty \prod\limits _{n = 1} ^ m  M _{\mu} ( 1 - B ^{-n}), 
\end{equation}
where $M _{\mu} $ is the moment generating function 
of the distribution $\mu$. This observation is not used
anywhere in the paper.

\begin{lem} \label{discombobulate}
	The series in \eqref{obtuse}
	converges for every $B>1$
	if and only if 
	for every $A >1$
	\begin{equation} \label{icky}
	\sum\limits _{m=1} ^\infty \prod\limits _{n = 1} ^ m 
	\sum \limits _{k = 0} ^{\infty} \mu (k) \left[ 1 - A ^{-n}\right]^{k}
	< \infty.
	\end{equation}
\end{lem}

\begin{proof}
	Suppose \eqref{obtuse} converges for every $B>0$. 
	It is sufficient to show that 
	\eqref{icky} holds for large $A$.
	Thus we can assume that $\mu([0,A]) \geq \frac 12$ and $A\geq {2}$, 
	so that $A ^n \geq An$ for $n \in \N$.
	We have 
	\begin{multline}
	\sum \limits _{k = 0} ^{\infty} \mu (k) \left[ 1 - A ^{-n}\right]^{k} 
	\leq  
	\sum \limits _{k = 0} ^{\infty} \mu (k) e^{- k A^{-n}}
	=  \sum \limits _{k = 0} ^{ \left\lceil A ^{2n} \right\rceil} \mu (k) e^{- k A^{-n}} +
	\sum \limits _{k = \left\lceil A ^{2n} \right\rceil + 1} ^{ \infty}  \mu (k) e^{- k A^{-n}}
	\\ 
	\leq \mu([0, A ^{2n}]) +  e^{- A^{n}}.
	\end{multline}
	For $m \in \N$, 
	\begin{equation}
	\frac{ \prod\limits _{n = 1} ^ m \left(\mu([0, A ^{2n}]) +  e^{- A^{n}} \right) }
	{ \prod\limits _{n = 1} ^ m \mu([0, A ^{2n}])   }
	\leq \prod\limits _{n = 1} ^ m \left(1 +   2 e^{- A^{n}} \right)
	< \prod\limits _{n = 1} ^ \infty \left(1 +   2 e^{- A^{n}} \right) < \infty.
	\end{equation}
	Hence 
	\begin{equation}
	\sum\limits _{m=1}  \prod\limits _{n = 1} ^ m \left(\mu([0, A ^{2n}]) +  e^{- A^{n}} \right)
	\simeq \sum\limits _{m=1}  \prod\limits _{n = 1} ^ m \mu([0, A ^{2n}]) < \infty,
	\end{equation}
	and 
	\begin{equation}
	\sum\limits _{m=1} ^\infty \prod\limits _{n = 1} ^ m 
	\sum \limits _{k = 0} ^{\infty} \mu (k) \left[ 1 - A ^{-n}\right]^{k}
	\leq 
	\sum\limits _{m=1}  \prod\limits _{n = 1} ^ m \left(\mu([0, A ^{2n}]) +  e^{- A^{n}} \right)
	< \infty.
	\end{equation}
	
	Conversely, 
	suppose \eqref{icky} holds for every $A>1$.
	Then for $B > 1$
	\begin{multline}
	\infty > 
	\sum\limits _{m=1} ^\infty \prod\limits _{n = 1} ^ m 
	\sum \limits _{k = 0} ^{\infty} \mu (k) \left[ 1 - 2 ^{-n} B ^{-n}\right]^{k}
	\geq 
	\sum\limits _{m=1} ^\infty \prod\limits _{n = 1} ^ m 
	\sum \limits _{k = 0} ^{\left\lfloor B ^{n} \right\rfloor} \mu (k) \left[ 1 - 2 ^{-n} B ^{-n}\right]^{k}
	\\
	\geq 
	\sum\limits _{m=1} ^\infty \prod\limits _{n = 1} ^ m 
	\sum \limits _{k = 0} ^{\left\lfloor B ^{n} \right\rfloor} \mu (k) \left[ 1 - 2 ^{-n} B ^{-n} k\right]
	\geq 
	\sum\limits _{m=1} ^\infty \prod\limits _{n = 1} ^ m 
	\mu ([0, B^{n}]) \left[ 1 - 2 ^{-n} \right]
	\\
	\geq
	\left( \prod\limits _{n = 1} ^ \infty \left( 1 - 2 ^{-n} \right) \right) ^{-1}
	\sum\limits _{m=1} ^\infty \prod\limits _{n = 1} ^ m 
	\mu ([0, B^{n}]).
	\end{multline}
\end{proof}

The next lemma helps in translating conditions 
\eqref{charade} and \eqref{lapidation}
of Propositions \ref{prop linear growth}
and \ref{prop superlinear growth}, respectively,
into conditions on $\mu$.
It is relevant 
for inequality \eqref{ply}
that 
for every $\varepsilon \in (0,1)$, $\varepsilon e^{1-\varepsilon} < 1$.

\begin{lem} \label{guff}
	For $n \in \N $ and $\varepsilon \in (0,1)$ 
	
	\begin{equation}\label{ply}
	\frac{\varepsilon ^ n e^{(1-\varepsilon)n}}{e 2^n \sqrt{ n}}
	\leq \P \Big\{\exists t \geq 0: \frac{S_t}{t} \geq \varepsilon ^{-1}, S _t \geq n \Big\}
	\leq 
	\frac{1}{(1 - \varepsilon)(1 -  \varepsilon  e^{1- \varepsilon})  }
	\frac{\big[ \varepsilon  e^{1- \varepsilon} \big] ^{n}}
	{ \sqrt{n} }.
	\end{equation}
\end{lem}
\begin{proof}
	On the one hand, 
	\begin{multline}
	\P \Big\{\exists t \geq 0: \frac{S_t}{t} \geq \varepsilon ^{-1}, S _t \geq n \Big\} 
	\geq 
	\PP{ \frac{\tau _n}{n} \leq \varepsilon  }  \PP{ S_{\tau _j} - S_{\tau _j - }  = 1, 
		j = \overline{1, n} }
	\\
	\geq e ^{-n \varepsilon} \frac{ \varepsilon ^n n ^n}{n! }  2 ^ {-n}
	\geq e ^{-n \varepsilon} \frac{ \varepsilon ^n n ^n e ^n}{ 2^n e n ^n  \sqrt{ n} }
	= \frac{ \varepsilon ^n  e ^{ (1-\varepsilon)n}}{e 2^n\sqrt{ n} }.
	\end{multline}
	\begingroup
	\allowdisplaybreaks
	On the other hand
	since 
	the sum of independent exponentials 
	have Erlang distribution
	for $m \in \N$
	\begin{align*}
	\PP{ \frac{\tau _{m}}{m} \leq \varepsilon   } &
	= \sum\limits _{i = m} ^{\infty} \frac{1}{i!}e^{- m\varepsilon}(m \varepsilon)^i
	\\ 
	&
	\leq e^{- m\varepsilon}  \sum\limits _{i = m}\frac{1}{m! m^{i-m}}(m \varepsilon)^i
	\\
	&
	= e^{- m\varepsilon} (m \varepsilon) ^m \frac{1}{m!}
	\sum\limits _{i = 0}\frac{1}{ m^{i}}(m \varepsilon)^i
	\\
	&
	\leq e^{- m\varepsilon} m^m \varepsilon ^m \frac{e^m}{m ^m \sqrt{m}} \frac{1}{1 - \varepsilon}
	=\frac{1}{(1 - \varepsilon)\sqrt{m}} \big[ \varepsilon  e^{1- \varepsilon} \big] ^{m}
	\end{align*}   
	and hence
	\begin{align*} 
	\P \Big\{\exists t \geq 0: \frac{S_t}{t} \geq \varepsilon ^{-1}, S _t \geq n \Big\} 
	= & \P \Big\{\exists k \geq 0: \frac{S_{\tau _{n+k}}}{\tau _{n+k}} \geq \varepsilon ^{-1}, S _{\tau _{n+k}} \geq n \Big\} \notag
	\\
	\leq & 
	\P \Big\{\exists k \geq 0: \frac{S_{\tau _{n+k}}}{\tau _{n+k}} \geq \varepsilon ^{-1} \Big\}
	\notag
	\\
	\leq & \sum\limits _{k = 0} ^\infty 
	\P \Big\{ \frac{\tau _{n + k}}{n+k} \leq \varepsilon  \Big\} 
	\notag 
	\\
	\leq  &
	\sum\limits _{k = 0} ^\infty 
	\frac{1}{(1 - \varepsilon) \sqrt{n+k} }
	\big[ \varepsilon  e^{1- \varepsilon} \big] ^{n+k}
	\notag
	\\
	\leq & 
	\frac{1}{(1 - \varepsilon) \sqrt{n} }
	\frac{\big[ \varepsilon  e^{1- \varepsilon} \big] ^{n}}
	{1 -  \varepsilon  e^{1- \varepsilon} }. 
	\end{align*}
\end{proof}
\endgroup

The next lemma gives a direct link between $\mu$
and the distribution of $\ell^{(A)} _n$. 
\begin{lem}\label{pungent}
	For $A>1$ 
	there exist $D_1 = D_1 (A)$, $D_2 =  D_2(A)$
	such that 
	\begin{equation}\label{clumpy}
	1 -    \sum\limits _{k = 0} ^ \infty \mu (k) \left[1 - D_1 ^{-n} \right] ^{k}
	\leq	r^{(A)} _n \leq 
	1 -   \sum\limits _{k = 0} ^ \infty \mu (k) \left[1 - D_2 ^{-n} \right] ^{k}, 
	\ \ \ n \in \N,
	\end{equation}	
	and $D_1(A), D_2(A) \to \infty$, $A \to \infty$. 
\end{lem}
\begin{proof}
	By 
	the definition of $\ell ^{(A)} _x$
	\begin{equation*}
	r^{(A)} _n =\PP{ \ell ^{(A)} _x \geq n  } = 1 - 
	\PP{ \ell ^{(A)} _x < n  }
	= 1 - \sum\limits _{k = 0 } ^\infty \mu (k)\left[
	\P \Big\{\forall t \geq 0: \frac{S_t}{t} < A \text{ or }
	S_t < n  \Big\}\right] ^k.
	\end{equation*}
	Note that $\left\{\forall t \geq \tau _n, 
	\frac{S_t}{t} < A \text{ or }
	S_t < n  \right\}$
	is the complement of the event 
	on the left hand side of \eqref{ply}
	with $A = \varepsilon ^{-1}$.
	Thus \eqref{clumpy} follows from Lemma \ref{guff}.
\end{proof}

\begin{lem}\label{shoehorn =? pigeonhole}
	Let $A > 1$.
	The convergence of the infinite product
	\begin{equation}
	\prod\limits_{  k =1 } ^\infty (1 - r_k ^{(A)}) > 0
	\end{equation}
	is equivalent to $	\sum\limits _{k = 1}  \mu(k) \log k  < \infty $.
	Likewise, 
	\begin{equation}
	\prod\limits_{  k =1 } ^\infty (1 - r_k ^{(A)}) = 0
	\end{equation}
	if and only if 
	$	\sum\limits _{k = 1}  \mu(k) \log k = \infty $.
\end{lem}
\begin{proof}
	It follows from Lemma \ref{pungent} 
	that \eqref{sulky} holds with $r_n = r_n ^{(A)}$,
	and hence 
	the statement follows from 
	Lemma \ref{anodyne}.
\end{proof}

\begin{lem}\label{invigilate=supervise drng examination}
	The convergence \eqref{lapidation} takes
	place for all $A>0$
	if and only if 
	\eqref{obtuse}
	holds for all $B>1$.
\end{lem}
\begin{proof}
	Suppose \eqref{lapidation} takes
	place for all $A>0$.
	By Lemma \ref{pungent}
	for $B >1$
	there exists 
	$A > 1$ 
	such that 
	\begin{equation}
	1 - r_n ^{(A)} \geq \sum\limits _{k = 0} ^\infty \mu (k) \left[ 1 - B ^{-n}  \right].
	\end{equation}
	Hence
	\begin{equation} 
	\infty > \sum\limits_{  m =1 } ^\infty \prod\limits_{  n =1 } ^m (1 - r_n ^{(A)}) 
	\geq \sum\limits_{  m =1 } ^\infty  \prod\limits _{n=1} ^ m 
	\sum\limits _{k = 0} ^ \infty \mu (k) \left[1 - B ^{-n} \right].
	\end{equation}
	Since $B > 1$ is arbitrary, the convergence  \eqref{obtuse} for all $B >1$ follows
	from Lemma \ref{discombobulate}.
	
	
	Conversely,
	suppose \eqref{obtuse}
	holds  for all $B>1$.
	By Lemma \ref{pungent} for some $A_1 >1$,
	\begin{equation} 
	\sum\limits_{  m =1 } ^\infty \prod\limits_{  n =1 } ^m (1 - r_n ^{(A)}) \leq 
	\sum\limits_{  m =1 } ^\infty  \prod\limits _{n=1} ^ m 
	\sum\limits _{k = 0} ^ \infty \mu (k) \left[1 - A_1 ^{-n} \right].
	\end{equation}
	The latter series converges by Lemma \ref{discombobulate}, 
	hence  \eqref{lapidation}
	holds for  $A >1$.
\end{proof}

We now come to the final part of the section.

\begin{proof}[Proof of Theorem \ref{thm linear growth d = 1}]
	Due to symmetry
	it suffices to show 
	that
	the   spread is linear in direction $+\infty$ only. 
	Indeed, by considering the simple 
	random walks $\{ ( - S _{t} ^{(x ,j)}, t \geq 0 ), x \in \Z ^\d, j \in \N \}$
	instead of $\{ (S _{t} ^{(x ,j)}, t \geq 0 ), x \in \Z ^\d, j \in \N \}$,
	we get another one-dimensional stochastic combustion growth process
	which is a reflection of the original one with respect to the origin.
	If the reflected process spreads linearly in direction $+\infty$,
	the original process spreads linearly in direction $-\infty$.
	The statement of the theorem thus follows from Proposition \ref{prop linear growth}
	and 
	and Lemma \ref{shoehorn =? pigeonhole}. 
\end{proof}

\begin{rmk}
	We see from the proof of Proposition \ref{prop linear growth}
	that in the settings of Theorem \ref{thm linear growth d = 1}  the growth toward $+\infty$
	would still remain linear
	even if all particles left of the origin were activated at time $t=0$.
	Of course, to consider
	such an initial configuration one would  need to 
	construct rigorously the process started with infinitely many active particles.
	If
	a.s. for all $t \geq 0$,
	$$ \sup\{ x + S _{s} ^{(x ,j)}: x \leq 0, j = 1,...,\eta(x), 0\leq s \leq t  \} < \infty,$$
	the construction may follow the standard arguments as a.s. only finitely
	many new sites will be activated within finite time intervals.
\end{rmk}

To complete the proof of  (ii) 
of
Theorem
\ref{thm main} we need Proposition \ref{prop superlinear growth higher dimension}.

\begin{proof}[Proof of Proposition \ref{prop superlinear growth higher dimension}]
	\label{page proof prop superlinear growth higher dimension}
	The projections of
	active particles in a
	$\d$-dimensional 
	stochastic combustion growth process
	on the first coordinate axis 
	perform a slowed down 
	simple
	continuous-time  random walk. 
	
	Given the set up for the $\d$-dimensional
	stochastic combustion growth processes,
	we are going to construct 
	a slowed down copy 
	of the 
	one-dimensional
	stochastic combustion growth processes
	which spreads slower than the  $\d$-dimensional
	process. 
	For $n \in \Z_+$,
	let $ \gamma _n = \min\{t \geq 0:
	\A _t \cap \{ (x_1, ..., x_\d) \in \Z^\d:
	x_1 = n \} \}$, that is,
	$ \gamma _n$
	is the  moment when an active particles 
	enters the plane $P_n :=\{ (x_1, ..., x_\d) \in \Z^\d:
	x_1 = n \} $. Let $L_n$ be the location
	where $P_n$ is first visited by an active particle 
	(this happens at time $\gamma _n$). 
	Modify the $\d$-dimensional 
	stochastic combustion growth processes
	so that in each moment $\gamma_n$, $n \in \Z_+$,
	all sleeping particles are removed from 
	$P_n$, while the active particles remain untouched;
	$1 + \eta (L_n)$
	particles on $L_n$ at time $\gamma _n$ 
	are active and thus  retained.
	In particular, at time $\gamma _0$
	all sleeping particles
	are removed from $P_0$, thus every site
	of $P_0$ excluding the origin is vacated.
	The evolution then proceeds as follows:
	one of the active particles started in the origin
	eventually hits either $P_{-1}$ or $P_1$;
	let us say 
	$P_1$ is hit at location $L_1 = (1, x_2,...,x_\d)$.
	The sleeping particles at  $L_1$ are activated
	at time $\gamma_1$, the sleeping 
	particles on $P_1\setminus \{L_1\}$
	are removed at the same time, and so on.
	
	In the resulting `trimmed'
	process the projections 
	of the active particles
	perform a  one-dimensional 
	stochastic combustion growth processes
	with particles jumping at  rate  $\frac{1}{\d}$
	instead of $1$. 
	Since the trimmed process
	is slower than the original $\d$-dimensional
	one,  the a.s. superlinear
	spread for
	a one-dimensional  stochastic combustion growth processes
	implies 
	the
	superlinear
	spread 
	for the $\d$-dimensional
	stochastic combustion growth processes.
\end{proof}

\begin{proof}[Proof of (ii) 
	of
	Theorem
	\ref{thm main}]
	The theorem is then a consequence of Propositions
	\ref{prop superlinear growth higher dimension} and 
	\ref{prop superlinear growth}
	and Lemma \ref{invigilate=supervise drng examination}.
\end{proof}


\section{Proof of  (iii)
	of
	Theorem
	\ref{thm main}} \label{sec proofs TADBP indirect}

By Proposition \ref{prop superlinear growth higher dimension}
it is enough to consider the case $\d=1$ only.
In this section  we give the proof
of  (iii)
of
Theorem
\ref{thm main}
for the one-dimensional system.

Recall that  $\sigma _x = \min\{t \geq 0: x \in \A _t  \}$
is the moment when $x$ is visited by an active particle
for the first time.
For $a > 0$ 
denote by $\chi _a$  the first time
when a simple continuous-time random walk started at $0$
hits $[a, \infty)$. Set $u =  \frac{1}{\sqrt{2\pi }} \int\limits_{-1} ^1 e^{-\frac{t^2}{2}} dt$
and let $ \varepsilon \in (0,1)$ be a small constant. By the central limit theorem
and reflection principle, for large $q$
\begin{equation} \label{scrim fabric}
u(1-\varepsilon)\leq	 \PP{ \chi _{q^{1/2}} \geq q } \leq u(1 + \varepsilon).
\end{equation}

Note that for $q > 1$, $x \geq q ^{1/2}$, $y \in [x - q^{1/2}, x-1]$
\begin{equation}\label{limber}
\{ \sigma_{x} - \sigma _{x - 1} \geq q  \}
\subset \{ S^{y, j} _{q} \leq x - y \text{ for all } j = 1,\dots \eta (y)  \}.
\end{equation}

Consequently for sufficiently large $q$
for $x \geq q^{1/2}$

\begin{equation}\label{walk out on updated}
\begin{split}
\PP{ \sigma_{x} - \sigma _{x - 1} \geq q }
& \leq \PP{ S^{y, j} _{q} \leq x - y \text{ for all } y = x-1, \dots, x - \lfloor q^{1/2} \rfloor,
	j = 1,\dots \eta (y)  }
\\
&  = \prod\limits _{y = x - \lfloor q^{1/2} \rfloor} ^{x-1} 
\sum\limits _{k = 0} ^\infty \mu (k) \left[ \PP{ \chi _{x - y} \geq q  } \right] ^k
= \prod\limits _{z =1   } ^{\lfloor q^{1/2} \rfloor} 
\sum\limits _{k = 0} ^\infty \mu (k) \left[ \PP{ \chi _{z} \geq q  } \right] ^k
\\
&
\leq  
\prod\limits _{z =1   } ^{\lfloor q^{1/2} \rfloor} 
\sum\limits _{k = 0} ^\infty \mu (k) \left[ \PP{ \chi _{q^{1/2}} \geq q  } \right] ^k
\leq
\bigg(  \sum\limits _{k = 0} ^\infty \mu (k) (1 + \varepsilon) ^k u ^k   \bigg) ^{q^{1/2} }
\end{split}
\end{equation}

By
\eqref{walk out on updated}
for sufficiently large $q$
\begin{equation} \label{fungible}
\PP{ \sigma_{x} - \sigma _{x - 1} \geq q }
\leq  c_\varepsilon ^{q^{1/2}} = \exp\left\{ - |\ln  c_ \varepsilon| q^{1/2}   \right\}.
\end{equation}
where $c_ \varepsilon  = \sum\limits _{k = 0} ^\infty \mu (k) (1 + \varepsilon) ^k u ^k   < 1$
for sufficiently small $\varepsilon$.

By \eqref{fungible}
for $c>2$
\begin{multline} \label{connoisseur}
\sum\limits_{ x \in \N } \P \Big\{ \max\limits_{1 \leq y \leq x} (\sigma_{y} - \sigma _{y - 1}) \geq  c^2 \left(\frac{\ln  x}{\ln  c_ \varepsilon}\right) ^2 \Big\}
\leq
\sum\limits_{ x \in \N } 
\sum\limits_{y = 1} ^x \P \Big\{  \sigma_{y} - \sigma _{y - 1} \geq  c^2 \left(\frac{\ln  x}{\ln  c_ \varepsilon}\right) ^2 \Big\} 
\\
\precsim
\sum\limits_{ x \in \N } x \exp\Big\{ - |\ln  c_ \varepsilon| c  \frac{\ln  x}{|\ln  c_ \varepsilon|}  \Big\}
= \sum\limits_{ x \in \N } x ^{-(c-1)} < \infty.
\end{multline}
Hence for $c>2$
\begin{equation}\label{morose}
\P \bigg\{   \max\limits_{1 \leq y \leq x} (\sigma_{y} - \sigma _{y - 1}) \geq  c^2 \left(\frac{\ln  x}{\ln  c_ \varepsilon}\right) ^2 \text{ infinitely often} \bigg\} = 0.
\end{equation}
Later  
\eqref{morose}
will be used 
for bounding the time it takes to travel
across `slow' sites.

Let $A > 1$.
Consider  
TADBP
on $\Z _+$
with $\psi _x = \ell ^{(A)} _x$.
It is possible that there exists an unbounded component,
that is, there is $y \in \N$ such that $y \xrightarrow{\Z _+} \infty$.
In this case we may proceed
as in the proof of Proposition \ref{prop superlinear growth}.
In the rest of the present proof
we  exclude this
case and assume
that
a.s.
no site is connected to infinity.
Note  that
under 
assumptions in (iii)
of
Theorem
\ref{thm main},
the existence of
$y \in \Z _+$
satisfying
$y \xrightarrow{ \Z_+} + \infty $ is not guaranteed;
see Proposition \ref{independent series}
and characterization of connected components in  TADBP
in Section \ref{sec TADBP percolation}.

Denote by $R_n$ the rightmost site of the $n$-th  connected component
(counting from the origin to the right) in the realization of 
TADBP
on $\Z _+$
with $\psi _x = \ell ^{(A)} _x$. 
Recall that 
$\ell ^{(A)} _x$ is defined in 
\eqref{secrete early}.
Note that $R_n$ depends on $A$. 
The number of site
in the interval $[0, R_n]$
which are dry is $n-1$ (specifically, the dry sites
are the leftmost sites of every component starting from the second;
recall that we consider the origin to be wet).
Denote by $l_k$ the length of $k$-th connected component,
$ l_k = R_k - R_{k-1} - 1 $.
The  random variables $\{l_k\}_{k \in \N}$
are
the excursions from $0$ of the Markov chain $\{Y_m\}_{m \in \Z_+}$
defined in \eqref{commandment}.
In particular, a.s. $l_1 \geq  \ell ^{(A)}  _0 $,
since the length of the first component 
is at least $\ell ^{(A)} _0$. 

By  Lemma \ref{pungent} for some $D = D(A)$
\begin{align*}
\sum\limits _{n = 1} ^ \infty \PP{l_n > n \ln ^2 n }
&
\geq 
\sum\limits _{n = 1} ^ \infty \PP{\ell ^{(A)} _n > n \ln ^2 n }
= 
\sum\limits _{n = 1} ^ \infty r ^{(A)} _{ \lceil  n \ln ^2 n \rceil}
\\
&
\geq \sum\limits _{n = 1} ^ \infty 
\sum\limits _{k = 0} ^ \infty  \mu (k) \left( 1 - \left[1 - D ^{-\lceil  n \ln ^2 n \rceil} \right] ^{k} \right)
\notag
\\
&
\geq
\sum\limits _{n = 1} ^ \infty 
\sum\limits _{k = 0} ^ \infty  \mu (k) \left( 1 - e^{-k D ^{-\lceil  n \ln ^2 n \rceil}}  \right)
\\
&
\geq
\sum\limits _{n = 1} ^ \infty  (1 - e^{-1}) \mu([ D ^{\lceil  n \ln ^2 n \rceil}, \infty ) ).
\end{align*}
Hence
by \eqref{enclave}
\begin{equation}
\sum\limits _{n = 1} ^ \infty \PP{l_n > n \ln ^2 n } = \infty,
\end{equation}
and consequently (\cite[Lemma 1]{CR61})
\begin{equation}
\PP{ \limsup _{n \to \infty } \frac{ l_n}{ n \ln ^2 n}  = \infty } = 1.
\end{equation}
In particular, a.s. $\limsup _{n \to \infty} \frac{R _n}{n \ln ^2 n} = \infty$.
Since the map 
$$ (10, \infty) \ni x \mapsto \frac{x}{\ln ^2 x} \in \R _+$$
is an increasing function,
for any $C > 1$
\begin{equation}
\limsup _{n \to \infty} \frac{R _n}{ n \ln ^2 R _n ^2} \geq 
\liminf _{n \to \infty} \frac{C n \ln ^2 n}{ n \ln ^2 \left(C ^2 n ^2 \ln ^4 n  \right) }
=
\liminf _{n \to \infty} \frac{C  \ln ^2 n}{   \left( 2 \ln C +  2\ln  n + 4 \ln  \ln n  \right) ^2 } = \frac{C}{4}.
\end{equation}
That is, 
\begin{equation}\label{harrowing}
\limsup _{n \to \infty} \frac{R _n}{ n \ln ^2 R _n ^2}  = \infty.
\end{equation}
By \eqref{morose} for large $n$ for all $y \in [0, R_n]$

\begin{equation}
\sigma_{y} - \sigma _{y - 1} \leq  c_2^2 \left(\frac{\ln  R_n}{\ln  c_ \varepsilon}\right) ^2,
\end{equation}
where $c_2 > 2$.
Using the same arguments
as when Lemma \ref{hold feet to the fire}
was applied in the proof of 
Proposition \ref{prop superlinear growth}, 
by  Lemma 
\ref{hold feet to the fire2}
we see that
the $n$-th
connected component is traversed 
within at most $\frac{2 l_n}{A}$
units of time. 
Hence the time needed to reach $R_n$
\begin{equation}
\sigma _{_{R_n}} \leq \frac{2}{A} \sum\limits _{k = 1} ^ n l_k + (n-1) c_2^2 \left(\frac{\ln  R_n}{\ln  c_ \varepsilon}\right) ^2
< \frac{2 R_n}{A}  + n c_2^2 \left(\frac{\ln  R_n}{\ln  c_ \varepsilon}\right) ^2,
\end{equation}
and by \eqref{harrowing}
\begin{multline}
\limsup _{n \to \infty } \frac{R_n}{ \sigma _{R_n}}   \geq 
\limsup _{n \to \infty } \frac{R_n }{ \frac{2}{A} R_n + n c_2^2 \left(\frac{\ln  R_n}{\ln  c_ \varepsilon}\right) ^2}
=  \limsup _{n \to \infty }
\frac{1 }{ \frac{2}{A}  + \frac{n \ln ^2  R_n}{R_n} \left(\frac{ c_2}{\ln  c_ \varepsilon}\right) ^2}
\\
=
\frac{1 }{ \frac{2}{A}  +  \left(\frac{ c_2}{\ln  c_ \varepsilon}\right) ^2
	\liminf _{n \to \infty }\frac{n \ln ^2  R_n}{R_n} } = \frac A2.
\end{multline}
Since $A > 1$ is arbitrary, it follows that 
\begin{equation}
\limsup _{m \to \infty } \frac{m}{ \sigma _m} = \infty.
\end{equation}
\qed

\begin{rmk}\label{snorkel}
	Let $\d=1$.
	The case not covered by Theorem \ref{thm main}
	is when
	the series
	\begin{equation}\label{hermit}
	\sum\limits _{k = 1} ^\infty  \mu(k) \log k . 
	\end{equation}
	diverges but slowly, so that 
	\begin{equation*}
	\sum\limits _{n \in \N}   \mu \left( [ e ^ {n \ln ^2 n}, \infty) \right)  < \infty
	\end{equation*}
	and 
	for some 	 $B >1$
	\begin{equation*}  
	\sum\limits _{m=1} ^\infty \prod\limits _{n = 1} ^ m 
	\mu \left( [0, B ^ n] \right) = \infty.
	\end{equation*}
	
	One might be tempted to conjecture that
	the spread is superlinear 
	once \eqref{hermit} diverges.
	Indeed, as we saw earlier in the proofs
	(specifically, Lemma \ref{shoehorn =? pigeonhole} and Section \ref{sec TADBP percolation}),
	the following holds true:
	for every $A>1$
	the fraction of sites 
	that are traversed by a particle
	moving toward $+\infty$
	at speed at least $A$ 	(as defined in \eqref{secrete early})
	is one.
	It would therefore suffice to show that 
	the time needed to traverse the rest of the sites (the `slow' sites) is not too large.
	This is done in the proof of  (iii)
	of
	Theorem
	\ref{thm main}
	under an additional assumption \eqref{enclave}.
	A better bound on $ \max\limits_{1 \leq y \leq x} (\sigma_{y} - \sigma _{y - 1})$
	could allow to weaken \eqref{enclave}, 
	but it is unclear if \eqref{enclave} can be dispensed with altogether.
	It may be that if \eqref{hermit}
	diverges very very slowly, 
	the slow sites have enough of an effect to bog down the growth
	and the spread is linear; or it may be that the spread is superlinear
	no matter how slowly \eqref{hermit} diverges.
	Both possibilities  seem  
	plausible to the authors of this paper;
	if a guess (or conjecture) had to be made, the latter would be chosen. 
	We would also like to note that the slow sites
	should still be slow 
	if the initial configuration of particles is capped
	so that there are $\eta (x) \wedge M$ at $x$,  $M \in \N$.
	Thus, understanding the slow sites for the 
	stochastic combustion growth process
	with bounded initial configuration may prove 
	helpful in shedding light on the cases not covered by Theorem \ref{thm main}.

	The gap between the conditions for the linear spread and superlinear
	spread in Theorem \ref{thm main} widens as the dimension $\d$ increases.
	Indeed, our proofs of the superlinear spread
	rely on Proposition \ref{prop superlinear growth higher dimension}
	and thus 
	are essentially one-dimensional. 
	One might therefore hope
	to  weaken the conditions implying 
	the
	superlinear spread
	for $\d\geq 2$
	by using techniques that would take into account 
	the spatial structure of $\Z^\d$. 
	
\end{rmk}

\section{Proof of Theorem \ref{thm linear growth d geq 2}} \label{sec proofs animals}

Recall that $\d \geq 2$ in the settings of Theorem \ref{thm linear growth d geq 2}.
Let $\theta _1 ^{(x,j)} < \theta _2 ^{(x,j)} < \dots $
be the jump times of the random walk $ \{S _t ^{(x,j)}, t \geq 0 \}$.
Let  $A > 1$.
We will see later in the proof
that we need $A$ to be  large enough to satisfy \eqref{crest} below.
Define 
\begin{equation}
W _{_A} ^{(x,j)} = \max \Big\{ n \in \N:  \frac{\theta _n ^{(x,j)}}{n} \leq \frac 1A    \Big\}.
\end{equation} 
Since $A> 1$,  $W _{_A} ^{(x,j)}$ is a.s. finite. 
Note how  
$W _{_A} ^{(x,j)}$ defined here differs 
from $
\ell ^{(A)} _x
$  in \eqref{secrete early}:
the former is defined solely in terms
of the jump moments, the latter is not.
Define also
\begin{equation}
W _{_A} ^{(x)}   =  \max\limits _{1 \leq j \leq \eta (x)} W _{_A} ^{(x,j)}
= \max \Big\{  n \in \N:  \frac{\theta _n ^{(x,j)}}{n} \leq \frac 1A, j \in \{1, \dots, \eta (x)\}   \Big\}.
\end{equation} 
and 
\begin{equation}
\rho _{_A}   = \max \Big\{ n \in \N:  \frac{ \tau _n }{n} \leq \frac 1A    \Big\}.
\end{equation} 
Note that $\rho _{_A}$ is equal in distribution to  $W _{_A} ^{(x,j)}$.
For $r \in \N$
and $y \in \Z ^\d$
let $\mathring{B}_1(y, r)$
the closed $\ell _1$-ball in $\Z ^\d$ with radius $r$ around $y$
with the removed center:
\begin{equation}
\mathring{B}_1(y, r)= 
\{x \in \Z^\d: |x-y|_1 \leq r  \}
\setminus \{y\}.
\end{equation}
In this chapter the union
$\mathcal{U} := \bigcup\limits _{y \in \Z ^\d}  \mathring{B}_1(y,
W _{_A} ^{(y)})$ 
plays the role of `potentially fast sites':
it is not necessary 
that every site $x \in \mathcal{U}$
is traveled by a particle at an average speed at least $A$ 
in some direction,
but from the definition of 
$W _{_A} ^{(y)}$
it follows that 
a site
$\widetilde x \notin \mathcal{U}$
cannot be traveled at speed exceeding $A$,
in the sense that a.s. 
\[
\forall y \in \Z ^\d, y \ne \widetilde x,
i \in \{1,..., \eta (y)\}:
\text{ there is no } t \geq q \geq 0 \text{ satisfying }
\frac{J ^{(y,i)} _t}{t} \geq A \text{ and }  S ^{(y,i)} _q = 
\widetilde x,
\]	
where $J ^{(y,i)} _t$ is the number of jumps of the random walk 
$S ^{(y,i)} _t$ before $t$.

\begin{lem} \label{morass}
	Let $\varepsilon = \frac  1A \in (0,1)$. For large $n$
	\begin{equation}
	\P \{ \rho _{_A} \geq n  \} \leq (\varepsilon e^{1-\varepsilon}) ^n.
	\end{equation}
\end{lem}
\begin{proof}
	We  use the same arguments as for the second inequality in \eqref{ply}.
	For $n$ satisfying $(1-\varepsilon)(1-\varepsilon e^{1-\varepsilon})\sqrt{2 \pi n} > 1$
	we have
	\begin{multline*}
	\P \{ \rho _{_A} \geq n  \} = \P\Big\{ \exists m \geq n : \frac{m}{\tau _m} \geq A \Big\}
	= \PP{ \exists m \geq n : \tau _m \leq m\varepsilon }
	\leq \sum\limits _{m=n} ^\infty \PP{ \tau _m \leq m\varepsilon}
	\\
	=
	\sum\limits _{m=n} ^\infty e ^{-m \varepsilon} \sum\limits _{j = m} ^\infty  \frac{(m\varepsilon)^j}{j!} 
	\leq 
	\sum\limits _{m=n} ^\infty e ^{-m \varepsilon}
	\frac{(m\varepsilon)^m}{m!} 
	\sum\limits _{j = 0} ^\infty  \frac{(m\varepsilon)^j}{ m^j} 
	\leq
	\sum\limits _{m=n} ^\infty e ^{-m \varepsilon}
	\frac{(m\varepsilon)^m e ^m}{m^m \sqrt{2 \pi m}} 
	\times \frac{1}{1-\varepsilon}
	\\
	=  \frac{1}{1-\varepsilon}
	\sum\limits _{m=n} ^\infty e ^{-m \varepsilon}
	\frac{\varepsilon^m e ^m}{\sqrt{2 \pi m}} 
	\leq \frac{1}{1-\varepsilon} \times
	\frac{1}{\sqrt{2 \pi n}} \times
	\frac{\varepsilon^n e ^{n(1-\varepsilon)}}{ 1 - \varepsilon e ^{1-\varepsilon} }
	< (\varepsilon e^{1-\varepsilon}) ^n. 
	\end{multline*}
\end{proof}




The random variables $\{ W _{_A} ^{(x)} \}_{x \in \Z ^\d }$
are independent and identically distributed.
Let $W _{_A}$ be a copy of $W _{_A} ^{(x)}$,$x \in \Z ^\d $,
independent of the sequence $\{ W _{_A} ^{(x)} \}_{x \in \Z ^\d}$.
By Lemma \ref{morass}

\begin{multline*}
\P \{ W _{_A} \geq n  \} = 1 - \P \{ W _{_A} < n  \} = 
1 - \sum\limits _{k = 0} ^\infty \mu (k) \left( \P \{ \rho _{_A} < n  \}  \right) ^k
\\
= 
1 - \sum\limits _{k = 0} ^\infty \mu (k) \left( 1 - \P \{ \rho _{_A} \geq n  \}  \right) ^k
\leq 1 - \sum\limits _{k = 0} ^\infty \mu (k) 
\left( 1 - (\varepsilon e^{1-\varepsilon}) ^n \right) ^k
\\
\leq 
1 - \sum\limits _{k = 0} ^\infty \mu (k) 
\left(1 -  1 \wedge k(\varepsilon e^{1-\varepsilon}) ^n \right)
= \sum\limits _{k = 0} ^\infty \mu (k) 
\left( 1 \wedge   k(\varepsilon e^{1-\varepsilon}) ^n \right)
\end{multline*}
Letting $B = (\varepsilon e^{1-\varepsilon}) ^{-1} >1$
we get


\begin{align}
\sum\limits _{n \in \N} \left[\PP{ W _{_A} \geq n   } \right] ^{\frac 1 \d}
& \leq 
\sum\limits _{n \in \N} \bigg[\sum\limits _{k = 1} ^\infty \mu (k) 
\left( 1 \wedge   k  B ^{-n} \right)\bigg] ^{\frac 1 \d}  \label{day laborer}
\\ &
\leq  \sum\limits _{n \in \N} \Bigg[\sum\limits _{k = 1} ^{\lfloor B^{n} \rfloor} \mu (k) 
\left( 1 \wedge   kB ^{-n} \right)
+
\sum\limits _{k = \lceil B^{n} \rceil} ^{\infty} \mu (k) 
\left( 1 \wedge   kB ^{-n} \right)
\Bigg] ^{\frac 1 \d}
\notag
\\
&
= 
\sum\limits _{n \in \N} \Bigg[
B ^{-n}
\sum\limits _{k = 1} ^{\lfloor B^{n} \rfloor} k \mu (k) 
+
\mu(  [ B^{n} , \infty ) )
\Bigg] ^{\frac 1 \d} 
\notag
\\
&
\leq 2 ^{\frac  1 \d} 
\sum\limits _{n \in \N} \Bigg[
B ^{-n} 
\sum\limits _{k = 1} ^{\lfloor B^{n} \rfloor} k \mu (k) 
\Bigg] ^{\frac 1 \d} 
+
2 ^{\frac  1 \d}  \sum\limits _{n \in \N} \left[
\mu(  [ B^{n} , \infty ) )
\right] ^{\frac 1 \d}.
\notag
\end{align}
Denote $b _k = \mu ( ( B^{k-1} ,B^{k}  ] )$. We have
\begin{align*}
\sum\limits _{n \in \N} \Bigg[
B ^{-n} 
\sum\limits _{k = 1} ^{\lfloor B^{n} \rfloor} k \mu (k) 
\Bigg] ^{\frac 1 \d} 
&
\leq 
\sum\limits _{n \in \N} \left[
B ^{-n} 
\sum\limits _{m = 1} ^{n} b_m B ^ m 
\right] ^{\frac 1 \d} 
\\
&
\leq \sum\limits _{n \in \N}  \sum\limits _{m = 1} ^{n} b_m ^{\frac 1 \d}
B ^{-\frac n \d}  B ^{\frac m \d} 
=  \sum\limits _{m = 1} ^{\infty}
b_m ^{\frac 1 \d} B ^{\frac m \d} 
\sum\limits _{n = m} ^ \infty  B ^{-\frac n \d} 
\\
&
=
\sum\limits _{m = 1} ^{\infty}
b_m ^{\frac 1 \d} B ^{\frac m \d} 
\frac{B ^{-\frac m \d} }{ 1- B ^{-\frac 1 \d}} 
=
\frac{1 }{ 1- B ^{-\frac 1 \d}}  \sum\limits _{m = 1} ^{\infty}
b_m ^{\frac 1 \d}
\\
& 
\precsim
\sum\limits _{m \in \N} \left[
\mu(  [ B^{m} , \infty ) )
\right] ^{\frac 1 \d}.
\end{align*}
Therefore by \eqref{day laborer}

\begin{equation}
\sum\limits _{n \in \N} \left[\PP{ W _{_A} \geq n   } \right] ^{\frac 1 \d}
\precsim  \sum\limits _{m \in \N} \left[
\mu(  [ B^{m} , \infty ) )
\right] ^{\frac 1 \d},
\end{equation}
and
hence by \eqref{ad nauseam}
\begin{equation}
\sum\limits _{n \in \N} \left[\PP{ W _{_A} \geq n   } \right] ^{\frac 1 \d}
< \infty. 
\end{equation}
Since for each $n \in \N$, $\PP{ W _{_{A}} \geq n}  \xrightarrow{A \to \infty} 0$, 
by the monotone convergence theorem,  
\begin{equation}
\lim\limits _{A \to \infty}  \sum\limits _{n \in \N} \left[\PP{ W _{_{A}} \geq n   } \right] ^{\frac 1 \d}
= 0. 
\end{equation}
Combining this with Theorem 1.1 in  \cite{Martin02}
yields the existence of $A > 1$ such
that 
\begin{equation} \label{crest}
\limsup\limits _{n \to \infty} \sup_{x_0, \dots, x_n} \frac{1}{n+1}\sum\limits _{i = 0} ^{n}
W _{_A} ^{(x_i)} \leq \frac 13.
\end{equation}
where the supremum is taken over all
connected sets of $n+1$ elements of $\Z ^\d$
containing  the origin, that is,
$x_0 = \0, x_1, x_2, \dots, x_n \in \Z ^\d$, $x _i \ne x_j$, $i \ne j$,
and
$\min\limits _{j = 1,...,i-1}|x_{i} - x_j| = 1$, $i = 1,2,...,n$.


\begin{proof}[Proof of Theorem \ref{thm linear growth d geq 2}]
	Take an infinite sequence 
	$$\{(Z_n, t_n, i_n)\}_{n \in \Z _+}, \ \ \ (Z_n, t_n, i_n) \in \Z ^\d \times \R _+ \times \N, \ \ n \in N,$$
	$Z_0 = \0$, $t_0 = 0$,  such that 
	the particles at site $Z_{n+1}$
	are activated by the particle $(i_n, Z_n)$ that started at $Z_n$,
	and $t_n$ is the activation time for $Z_n$. 
	Let $z_0  = 0, z_1, z_2, ..., $ be the 
	successive sites 
	visited by the particles
	$(i_n, Z_n)$, $n \in \Z_+$,
	during the time interval $ [t_n, t_{n+1} ] $,
	so that
	$$ \bigcup\limits _{j = 0} ^{\infty}\{  z_j  \} = \bigcup\limits _{n = 0} ^{\infty} \{S^{(i_n, Z_n)}_t
	+ Z_n, \  0 \leq t \leq t_{n+1} - t_n \}$$
	and the sequence $\{  z_j  \}_{j \in \Z_+}$ 
	does not  contain repeating elements (that is,
	if a site is already in the sequence $\{  z_j  \}$, it is not appended
	even
	when  visited  by a particle  $(i_m, Z_m)$
	during  $ [t_m, t_{m+1} ]$).
	Note that  $ \{  Z_j  \}_{j \in \Z_+} \subset \{  z_j  \}_{j \in \Z_+}$. 
	
	Let $\text{trav}( Z_n, Z_{n+1}) $ 
	be the number of sites excluding $Z_n$
	visited by the particle $(i_n, Z_n)$
	by the moment  $t_{n+1}$ when it reaches $Z_{n+1}$  (for instance, if
	$Z_{n+1}$ and $Z_n$ are neighbors
	and 
	$(i_n, Z_n)$ goes directly from $Z_n$ to $Z_{n+1}$, then  $\text{trav}( Z_n, Z_{n+1}) = 1 $).
	For $n \in \N $ denote by 
	$\kappa_n \in \N$
	the index satisfying $Z_n = z_{\kappa_n}$.
	Note that $\kappa_n$ is uniquely defined. 
	As 
	the particle $(i_n, Z_n)$
	travels from $Z_n$ to $Z_{n+1}$,
	only those  sites are added
	to  $\{  z_j  \}_{j \in \Z_+}$
	which are not in the sequence already.
	Therefore a.s.
	\begin{equation*}
	\text{trav}( Z_n, Z_{n+1}) \geq \kappa_{n+1} - \kappa_n,
	\ \ \ n \in \Z_+
	\end{equation*} 
	and hence 
	\begin{equation}\label{workaday}
	\sum\limits _{k=0} ^{n-1} \text{trav}( Z_k, Z_{k+1}) 
	\geq  \kappa_n,  \ \ \ n \in \N.
	\end{equation} 
	By  \eqref{crest} a.s. for large $m \in \N$,
	\begin{equation}
	\frac{1}{m}\sum\limits _{j = 0} ^{m}
	W _{_A} ^{(z_j)} \leq \frac 12.
	\end{equation}
	Hence by \eqref{workaday} a.s. for large $n \in \N$
	\begin{equation}\label{mesh}
	\frac{\sum\limits _{j = 0} ^{n}
		W _{_A} ^{(Z_j)} }{\sum\limits _{j = 0} ^{n-1}
		\text{trav}( Z_j, Z_{j+1}) } \leq
	\frac{ \sum\limits _{j = 0} ^{\kappa_n}
		W _{_A} ^{(z_j)}} {\kappa_n }    \leq
	\frac 12.
	\end{equation}
	Next, note that if the path 
	from $Z_n$ to $Z_{n+1}$
	traveled by the particle $(i_n, Z_n)$
	is not entirely covered by 
	the set of `potentially fast' sites
	$\mathcal{U}$,
	then 
	\begin{equation}
	\text{trav}( Z_j, Z_{j+1}) \leq A (t_{j+1} - t_j).
	\end{equation}
	By \eqref{mesh}
	at least half of the sites  
	of the sequence $\{z_j\}_{j \in \Z_+}$
	are outside of $\mathcal{U}$, hence 
	a.s. for large $n$
	\begin{equation}\label{belabor}
	\sum\limits _{j = 0} ^{n-1} (t_{j+1} - t_j)
	\geq 
	\frac {1}{2A}
	\sum\limits _{j = 0} ^{n-1} \text{trav}( Z_j, Z_{j+1}).
	\end{equation}
	Note that by \eqref{crest},
	the inequality \eqref{belabor}
	holds uniformly in  $ \{  Z_j  \}_{j \in \Z_+}$; that is,
	there exists $n_0 \in \N$
	such that 
	for all choices of the sequence $Z_1, Z_2, ...$
	such that $Z_{n+1}$ activated by a particle
	started at $Z_n$,
	\eqref{belabor} holds for all  $n \geq n_0$.
	In other words, a.s.
	\begin{equation*}
	\liminf \limits _{n \to \infty}
	\inf\limits_{\{  Z_j  \}_{j \in \Z_+}}
	\frac{\sum\limits _{j = 0} ^{n-1} (t_{j+1} - t_j)}{\sum\limits _{j = 0} ^{n-1} \text{trav}( Z_j, Z_{j+1})} \geq \frac{1}{2A},
	\end{equation*}
	where the infimum is taken over all the sequences
	of successively activated sites as described 
	in the previous sentence. 
	Since $|Z_n|_1 \leq \sum\limits _{j = 0} ^{n-1} \text{trav}( Z_j, Z_{j+1})$, 
	we have
	by \eqref{belabor}
	for large $n$
	\begin{equation}
	\frac{t_n}{|Z_n|_1} \geq \frac{ \sum\limits _{j = 0} ^{n-1} (t_{j+1} - t_j) }{\sum\limits _{j = 0} ^{n-1} \text{trav}( Z_j, Z_{j+1})} 
	\geq \frac{ \frac 12 \sum\limits _{j = 0} ^{n-1}  \frac{\text{trav}( Z_j, Z_{j+1})}{A} }{\sum\limits _{j = 0} ^{n-1} \text{trav}( Z_j, Z_{j+1})} 
	=\frac{1}{2A},
	\end{equation}
	that is,
	a.s. for large $n $
	\begin{equation}\label{revamp}
	\frac{|Z_n|_1}{|t_n|} \leq 2A.
	\end{equation}
	Since a.s. \eqref{revamp}
	holds  uniformly for  every
	successively activated 
	sequence of sites, the proof is complete.
\end{proof}

\section{Convergence properties of related series}
\label{sec convergence properties}

In  this section we show that
the series in
$(ii)$ and $(iii)$
of Theorem \ref{thm main} have independent convergence properties. 
Thus, the content of this section
is not used in the proof  of Theorem \ref{thm main} per se
but rather addresses the logical independence of its parts.  
The construction in the proof of the following lemma is courtesy of Christian Remling.

\begin{lem} \label{tawdry}
	Let $\{ u_m\}_{m \in \N}$ be an increasing sequence of non-negative numbers, 
	$u_m \to \infty$, $m \to \infty$. Let $\{ v_m\}_{m \in \N}$
	be another sequence of positive numbers. There exists 
	a sequence $\{g _n\}_{n \in \mathbb{N}}$, $g _n \to 0$, $g _n > 0$ such that 
	\begin{equation} \label{hackneyed}
	\sum\limits_{n \in \mathbb{N} } g_n = \infty
	\end{equation}
	and 
	\begin{equation}\label{umbrage}
	\sum\limits _{m \in  \mathbb{N}} u_m \exp\Big\{ - \sum\limits _{i = 1} ^m g _i v_i \Big \} = \infty.
	\end{equation}
\end{lem}
\begin{proof}
	Define the elements of the sequences 
	$\{g _n\}_{n \in \mathbb{N}}$,
	$\{M _n\}_{n \in \mathbb{N}}$, and $\{K _n\}_{n \in \mathbb{N}}$
	consecutively
	as follows.
	Set  $K_1 = 1$ and $g_1 = 1$. For $n \in \N$ once $K_n$ and $g_i$, $i = 1,\dots,K_n$,
	are defined, set 
	\begin{equation} \label{a big ask}
	M_n = \min \bigg\{ m \geq K _{n} + 1: u_m \geq e \exp\Big\{  \sum\limits _{i = 1} ^{K_{n}} v_i g_i  \Big\} \bigg\}
	\end{equation}
	Then for $i = K_n + 1, \dots, M_n $, set $g_i = h_n:=  \Big((M_n - K_n) \max\limits _{ K_n +1 \leq j \leq M_n}  v_j \Big) ^{-1} \wedge \frac 1n $.
	A single step is completed by setting 
	$$
	K_{n+1} = \min\{ m \geq M_n +1: h_n(K_{n+1} - K_n  ) \geq 1 \},
	$$
	and $g_i = h_n$ for  $i = M_n + 1, \dots, K_{n+1} $.
	Next we define $M_{n+1}$ as in \eqref{a big ask} (of course
	with $n$ replaced by $n+1$ everywhere in \eqref{a big ask}),
	and so forth.

	With this construction we have 
	
	\begin{equation}
	\sum\limits_{i = K_n + 1 } ^{K_{n+1}} g_i \geq 1
	\end{equation}
	and 
	\begin{multline}
	u_{_{M_n}} \exp\Big\{ - \sum\limits _{i = 1} ^{M_n} g _i v_i \Big \} =
	u_{_{M_n}} \exp\Big\{ - \sum\limits _{i = 1} ^{K_n} g _i v_i \Big \}
	\times \exp\Big\{ - \sum\limits _{i = K_n + 1} ^{M_n} g _i v_i \Big \}
	\\
	\geq e \times \exp\Big\{ - h_n (M_n - K _n) \max\limits _{ K_n +1 \leq j \leq M_n}  v_j\Big \}
	\geq e e^{-1} = 1.
	\end{multline} 
	Thus both \eqref{hackneyed} and \eqref{umbrage} hold.
\end{proof}

\begin{prop} \label{independent series}
	The convergence properties of the series 
	in \eqref{obtuse}  and \eqref{enclave} are independent. 
	That is, all four combinations
	of both series converging,
	either one of the two converging,
	and both diverging are possible. 
\end{prop}
\begin{proof}
	Fix $B > 1$ and let $b _ n = \mu ([ 0, B^n])$
	and $c _n = \mu ([ 0, B^{n \ln  ^ 2 n} ))$.
	We only consider distributions $\mu$ with unbounded support.
	We have $b_n \nearrow 1$
	and
	\begin{equation}
	{b_{ \lfloor   n \ln  ^ 2 n  \rfloor } \leq c_n \leq b_{ \lceil   n \ln  ^ 2 n  \rceil }},
	\ \ \  n \geq 2.
	\end{equation}
	
	
	The series 	in \eqref{obtuse}  and \eqref{enclave}
	can be written as 
	$
	\sum\limits _{m=1} ^\infty \prod\limits _{n = 1} ^ m b_n  
	$
	and
	$
	\sum\limits _{n = 1} ^ \infty (1- c_n)
	$ respectively.
	Note that 
	\begin{equation} 
	\ln  ^ 2 k  \leq
	(k + 1) \ln  ^ 2 (k+1)  - k \ln  ^ 2 k 
	\leq \ln  ^ 2 k +  2 \ln k + 2, 
	\ \ \ k \in \N.
	\end{equation}
	Since the sequence $\{b_n\}$
	is monotone
	\begin{align}
	\sum\limits _{m=2} ^\infty \prod\limits _{n = 1} ^ m b_n 
	&	\geq 
	\sum\limits _{m=2} ^\infty \prod\limits _{k = 1} ^ { \max\{ l \in \N : l \ln ^2 l < m  \}} \prod\limits_{ 
		\substack{
			i \in \N : \\
			k\ln ^2 k \leq i < (k + 1) \ln  ^ 2 (k+1) 
		}
	} b_i \label{chuck}
	\\
	&	\geq
	\sum\limits _{m=2} ^\infty \prod\limits _{k = 1} ^ { \max\{ l \in \N : l \ln ^2 l < m  \}} \prod\limits_{ 
		\substack{
			i \in \N : \\
			k\ln ^2 k \leq i < (k + 1) \ln  ^ 2 (k+1) 
		}
	} c_k  \notag
	\\
	&	\geq 
	\sum\limits _{m=2} ^\infty \prod\limits _{k = 1} ^ { \max\{ l \in \N : l \ln ^2 l < m  \}}
	c_k ^{\ln  ^ 2 k +  2 \ln k + 2}  \notag
	\\
	&	\geq \sum\limits _{l=2} ^\infty 
	\ln ^2 l 
	\prod\limits _{k = 1} ^ { m}
	c_k ^{\ln  ^ 2 k +  2 \ln k + 2}   \notag
	\\
	&	=
	\sum\limits _{l=2} 
	\ln ^2 l  \exp\Big\{ -\sum\limits_{i=1} ^l \gamma _i (\ln  ^ 2 i +  2 \ln i + 2)    \Big\}, \notag
	\end{align}
	where $ \gamma _i : = - \ln c_i > 0$. Note that $\sum\limits _{n = 1} ^ \infty (1- c_n) \simeq \sum\limits _{n = 1} ^ \infty \gamma _n $
	since  $ \lim\limits _{n \to \infty} \frac{1 - c_n}{ \gamma _n} =1 $. 
	By Lemma \ref{tawdry}, $\{\gamma _i\} _{i \in \N}$ can be chosen in such a way that 
	
	\begin{equation*}
	\sum\limits _{m=1} ^\infty \prod\limits _{n = 1} ^ m b_n  = \infty
	\   \  \text{ and }  \  \  \  \sum\limits _{n = 1} ^ \infty (1- c_n)  = \infty.
	\end{equation*}
	
	The other three cases are more straightforward.
	Before proceeding to them, note that 
	similarly to \eqref{chuck}
	\begin{align}
	\sum\limits _{m=2} ^\infty \prod\limits _{n = 1} ^ m b_n 
	& \leq 
	\sum\limits _{m=2} ^\infty \prod\limits _{k = 1} ^ { \max\{ l \in \N : l \ln ^2 l < m  \} - 1} \prod\limits_{ 
		\substack{
			i \in \N : \\
			k\ln ^2 k  < i \leq(k + 1) \ln  ^ 2 (k+1) 
		}
	} b_i \label{notary}
	\\
	& \leq
	\sum\limits _{m=2} ^\infty \prod\limits _{k = 1} ^ { \max\{ l \in \N : l \ln ^2 l < m  \} - 1} \prod\limits_{ 
		\substack{
			i \in \N : \\
			k\ln ^2 k < i \leq (k + 1) \ln  ^ 2 (k+1) 
		}
	} c_{k+1} \notag
	\\ &
	\leq 
	\sum\limits _{m=2} ^\infty \prod\limits _{k = 1} ^ { \max\{ l \in \N : l \ln ^2 l < m  \} - 1}
	c_{k+1} ^{\ln ^2 k } \notag
	\\
	&
	\leq \sum\limits _{l=2} ^\infty 
	(\ln ^2 l +2 \ln l + 2 )
	\prod\limits _{k = 1} ^ { l}
	c_{k+1} ^{\ln ^2 k } \notag
	\\
	&
	=
	\sum\limits _{l=2} 
	(\ln ^2 l +2 \ln l + 2 )
	\exp\Big\{ -\sum\limits_{i=1} ^l \gamma _{i+1} \ln ^2 i         \Big\} \notag
	\end{align}
	Taking $\{\gamma _i\} _{i \in \N}$ very small (for example $\gamma _i = e ^{- i^2}$ ) we can easily achieve 	 
	\begin{equation*}
	\sum\limits _{m=1} ^\infty \prod\limits _{n = 1} ^ m b_n  = \infty
	\   \  \text{ and }  \  \  \  \sum\limits _{n = 1} ^ \infty (1- c_n)  < \infty.
	\end{equation*}
	Letting  $\{\gamma _i\} _{i \in \N}$ converge to $0$ very slowly,
	for example $\gamma _i = \frac{1}{\ln \ln \ln i} $ for large $i$,
	we get 
	\begin{equation*}
	\sum\limits _{m=1} ^\infty \prod\limits _{n = 1} ^ m b_n  < \infty
	\   \  \text{ and }  \  \  \  \sum\limits _{n = 1} ^ \infty (1- c_n) = \infty.
	\end{equation*}
	Letting $\gamma _{n+1} = \frac{1}{ n (\ln n) ^{3/2}}$, $n \geq 2$, 
	we get $\sum\limits _{n = 1} ^ \infty \gamma _n < \infty$
	and for large $m$
	$$
	\sum\limits_{i=1} ^l \gamma _{i+1} \ln ^2 i 
	=   \sum\limits_{i=1} ^l \frac {\ln ^2 i }{ i (\ln i) ^{3/2}}
	\geq \frac 12  \ln ^{3/2} l,
	$$
	and 	hence 
	\begin{equation*}
	\sum\limits _{l=2} ^\infty
	(\ln ^2 l +2 \ln l + 2 )
	\exp\Big\{ -\sum\limits_{i=1} ^l \gamma _{i+1} \ln ^2 i   \Big\}
	\precsim   \sum\limits _{l=2}  ^\infty
	\frac{\ln ^2 l +2 \ln l + 2}{l ^{ \frac 12 \sqrt{\ln l }}} < \infty. 
	\end{equation*}
	Taking  \eqref{notary} into account we see that both series of interest converge 
	\begin{equation*}
	\sum\limits _{m=1} ^\infty \prod\limits _{n = 1} ^ m b_n  < \infty
	\   \  \text{ and }  \  \  \  \sum\limits _{n = 1} ^ \infty (1- c_n)  < \infty.
	\end{equation*}
\end{proof}

\section*{Acknowledgements}

We would like to thank Christian Remling for 
the construction used in the proof
of Lemma \ref{tawdry}
and 
Martin Zerner 
for bringing to our attention the works \cite{Lamp70} and \cite{Kel06}
and thus  helping 
to place Section \ref{sec TADBP percolation} in the context of existing research.
We would also like to thank the anonymous referee
whose comments helped to improve the paper.


\bibliographystyle{alphaSinus}
\bibliography{bibliography}

\end{document}